\documentclass[11pt]{article}

\usepackage{a4,amssymb,amsmath,amsthm,color,graphicx}

\def\eins{\mbox{1\hskip-0.24em l}}
\def\T{^{\sf T}}
\newcommand{\C}{ {\mathbb C} }
\newcommand{\R}{ {\mathbb R} }
\newcommand{\N}{{\mathbb N}}
\newcommand{\MM}{{\mathbb M}}
\newcommand{\Uu}{{\mathbb U}}

\newcommand{\LL}{{\cal L}}
\newcommand{\PP}{{\cal P}}
\newcommand{\PO}{{\mathbb P}}

\newcommand{\diag}{\,\mbox{diag}}

\newcommand{\cc}{{\bf c}}
\def\ones{\mbox{\normalfont{1\hskip-0.24em l}}}
\newcommand{\gdw}{\ \iff\ }

\newtheorem{remark}{Remark}[section]
\newtheorem{theorem}{Theorem}[section]
\newtheorem{lemma}{Lemma}[section]

\definecolor{otherblue}{rgb}{0,0.3,0.6}

\begin{document}
\title{Discrete Adjoint Implicit Peer Methods\\ in Optimal Control}
\author{Jens Lang \\
{\small \it Technical University Darmstadt,
Department of Mathematics} \\
{\small \it Dolivostra{\ss}e 15, 64293 Darmstadt, Germany}\\
{\small lang@mathematik.tu-darmstadt.de} \\ \\
Bernhard A. Schmitt \\
{\small \it Philipps-Universit\"at Marburg,
Department of Mathematics,}\\
{\small \it Hans-Meerwein-Stra{\ss}e 6, 35043 Marburg, Germany} \\
{\small schmitt@mathematik.uni-marburg.de}}
\maketitle

\begin{abstract}
It is well known that in the first-discretize-then-optimize approach in the control of ordinary differential equations the adjoint method may converge under additional order conditions only.
For Peer two-step methods we derive such adjoint order conditions and pay special attention to the boundary steps.
For $s$-stage methods, we prove convergence of order $s$ for the state variables if the adjoint method satisfies the conditions for order $s-1$, at least.
We remove some bottlenecks at the boundaries encountered in an earlier paper of the first author et al.
[J. Comput. Appl. Math., 262:73--86, 2014]
and discuss the construction of 3-stage methods for the order pair (3,2) in detail including some matrix background for the combined forward and adjoint order conditions. The impact of nodes having equal differences is
highlighted. It turns out that the most attractive methods are related to BDF. Three 3-stage methods are constructed which show the expected orders in numerical tests.
\end{abstract}


\section{Introduction}
In this paper, we are interested in the numerical solution of
the following ODE-constrained nonlinear optimal control problem:
\begin{align}
\mbox{minimize } C\big(y(T)\big) \label{OCprob_objfunc} &\\
\mbox{subject to } y'(t) =& \,f\big(y(t),u(t)\big),\quad
u(t)\in U,\;t\in(0,T], \label{OCprob_ODE}\\
y(0) =& \,y_0, \label{OCprob_ODEinit}
\end{align}
where the state $y(t)\in\R^m$, the control $u(t)\in\R^d$,
$f: \R^m\times\R^d\mapsto\R^m$, the objective function $C: \R^m\mapsto\R$, and the set of admissible controls $U\subset\R^d$ is closed and convex. Introducing for any $u\in U$ the normal cone mapping
\begin{align}
\label{def_cone}
N_U(u) =&\, \{ w\in\R^d: w^T(v-u)\le 0 \mbox{ for all } v\in U\},
\end{align}
the first-order optimality conditions read \cite{Hager2000,Troutman1996}
\begin{align}
y'(t) =& \,f\big(y(t),u(t)\big),\quad t\in(0,T],\quad y(0)=y_0, \label{KKT_state}\\
p'(t) =& \,-\nabla_y f\big(y(t),u(t)\big)\T p,\quad t\in[0,T),
\quad p(T)=\nabla_y C\big(y(T)\big)\T, \label{KKT_adj}\\
& \,-\nabla_u f\big(y(t),u(t)\big)\T p \in N_U\big(u(t)\big),\quad t\in[0,T]. \label{KKT_ctr}
\end{align}
Under appropriate regularity conditions, there exists a local solution
$(y^\star,u^\star)$ of the optimal control problem
\eqref{OCprob_objfunc}-\eqref{OCprob_ODEinit} and a Lagrange multiplier
$p^\star$ such that the first-order optimality conditions
\eqref{KKT_state}-\eqref{KKT_ctr} are necessarily satisfied at $(y^\star,u^\star,p^\star)$. If, in addition, the Hamiltonian
$H(y,u,p):=p\T f(y,u)$ satisfies a coercivity assumption, then these
conditions are also sufficient \cite{Hager2000}. The control
uniqueness property introduced in \cite{Hager2000} yieds the existence
of a locally unique minimizer $u=u(y,p)$ of the Hamiltonian over all
$u\in U$, if $(y,p)$ is sufficiently close to $(y^\star,p^\star)$.
Substituting $u$ in terms of $(y,p)$ in \eqref{KKT_state}-\eqref{KKT_adj},
gives then the two-point boundary value problem
\begin{align}
\label{RWPy}
y'(t)=&\,g\big(y(t),p(t)\big),\quad y(0)=y_0,\\
\label{RWPp}
p'(t)=&\,\phi\big(y(t),p(t)\big),\quad p(T)=\nabla_y C\big(y(T)\big)\T,
\end{align}
with the source functions defined by
\begin{align}
g(y,p) :=& f\big(y,u(y,p)\big),
\quad \phi(y,p) := -\nabla_yf\big(y,u(y,p)\big)\T p.
\end{align}
This boundary value problem plays a key role in any consistency
and convergence analysis. In what follows, we will assume sufficient
smoothness of the optimal control problem, so that the
elimination of the control as described above can be always applied.
\par
Here we will follow the \textit{first-disretize-then-optimize} approach,
i.e., the ODE system (\ref{OCprob_ODE})-(\ref{OCprob_ODEinit}) is first
discretized by applying an $s$-stage implicit Peer two-step method. This
leads to a finite dimensional optimal control problem, for
which the first-order discrete optimality system can be derived and
solved by existing optimization solvers such as nonlinear Newton-type
algorithms. In spite of the large size of the resulting problems, the
flexibility of this approach naturally allows the incorporation of additional
constraints and bounds. Further advantages are the direct use of automatic
differentiation techniques and the computation of discrete adjoints, which
are consistent with the discrete optimal control problem. Symmetric approximations
of Hessian matrices can be easily derived and result in a computational speedup.
\par
To ensure optimal order of
convergence to the infinite dimensional optimality system, the discrete
adjoint equations should represent a consistent approximation of its
continuous counterpart (\ref{KKT_adj}) -- a property that is refered to
as \textit{adjoint consistency}. Adjoint-consistent one-step Runge-Kutta methods
were studied by Hager \cite{Hager2000}, Sandu \cite{Sandu2006},
and Pulova \cite{Pulova2008}. The special class of symplectic partitioned
Runge-Kutta methods and additional order conditions were investigated by
Bonnans and Laurent-Varin \cite{BonnansLaurentVarin2006} and in a different
setting earlier by Murua \cite{Murua1997}. Later on, symplectic properties
for implicit-explicit Runge-Kutta methods in the context of optimal control
were analyzed by Herty et al. \cite{HertyPareschiSteffensen2013}.
Lang and Verwer \cite{LangVerwer2013}
showed for third-order W-methods that they also have to fulfill additional
consistency conditions in order to make them valuable for optimal control. Reverse mode automatic differentiation on explicit Runge–Kutta methods is an alternative
approach to derive consistent discrete adjoints as shown by Walther \cite{Walther2007}.
However, one-step methods might suffer from serious order reduction, especially
when they are applied to very stiff problems or large-scale ODE systems obtained from semi-discretizations of PDE systems with general boundary conditions.
\par
The situation is more complex for multistep methods. Here, the discrete adjoint schemes
of linear multistep methods are in general not consistent or show a
significant decrease of the approximation order, see Sandu \cite{Sandu2008}
and Albi et al. \cite{AlbiHertyPareschi2019}. Backward differentiation formula (BDF)
and Peer methods \cite{BeckWeinerPodhaiskySchmitt2012,SchmittWeinerErdmann2005} which are particularly suitable for large-scale, nonlinear and stiff systems
of ODEs keep their high order in the interior of the time domain, but the adjoint initialization steps are usually inconsistent approximations \cite{BeigelMommerWirschingBock2014,SchroederLangWeiner2014} and the numerical approximation of missing starting values has to be done with care. These inherent difficulties have limited the application of multistep
methods for optimal control problems in a first-discretize-then-optimize
solution strategy.
In this paper we will propose a novel approach to
overcome these structural deficiencies for both Peer and BDF methods.
\par
The numerical results for Peer methods in the previous paper by Schr\"oder et al. \cite{SchroederLangWeiner2014} were quite disappointing since the adjoint solutions of the complete boundary value problem  did converge with order one only.
We overcome some bottlenecks by the following measures:
\begin{enumerate}
\item
Using redundant formulations since equivalent versions of the forward scheme are not equivalent in the adjoint scheme.
\item
Modified first and last time steps give additional degrees of freedom.
\item
Using a general approximation $y_h(T)=\sum_{j=1}^s w_jY_{Nj}$ at the end point $T$ allows for more general nodes with $c_s\not=1$.
\end{enumerate}
With respect to the second item, we remind that for time step schemes used only once, the order may be lower by one than the order of the overall scheme.

\par
The paper continues in Section~\ref{Svoradj} with the description of the Peer method and derivation of the adjoint schemes including the boundary conditions.
The order conditions of both schemes are derived in Section~\ref{SOrdng}.
Aiming at high-order convergence it will be shown in Section~\ref{SKonvg} that for an $s$-step method order $O(h^s)$ can be shown for the $y$-variable if the solution for the $p$-variable has $O(h^{s-1})$ convergence at least.
Accordingly, the construction of 3-stage methods is based on a thorough discussion of methods with the global order pair (3,2) for the solution $y$ and the adjoint $p$.
This is also motivated by the fact that the order pair (3,3) can not be satisfied in our present setting.
Since this discussion shows a certain preference for nodes with flip symmetry this question is pursued in Section~\ref{Sfao} in detail by combining the forward and adjoint order conditions.
Numerical tests in Section~\ref{STests} confirm the convergence results from Section~\ref{SKonvg}.
\section{Implicit Peer two-step methods: the forward and
adjoint schemes}\label{Svoradj}
Purely implicit Peer methods were introduced in \cite{SchmittWeinerErdmann2005} in a special form suited for parallel implementation.
They have a two-step structure on a time grid $\{t_0,\ldots,t_{N+1}\}\subset[0,T]$ with step sizes $h_n=t_{n+1}-t_n$ and use $s$ solution approximations $Y_{ni}\approx y(t_n+c_ih_n),\,i=1,\ldots,s$, per time step associated with a set of fixed off-step nodes $c_1,\ldots,c_s$.
The form of this scheme is not unique, an equivalent formulation saving the memory for function evaluations $F_{ni}:=f(Y_{ni})$ was used in other papers with the first author, e.g. \cite{SchroederLangWeiner2014}.
Both versions of the scheme produce identical approximations $Y_{ni}$ and one may choose either one for the integration forward in time.
However, this is no longer true if adjoint equations come into play.
Hence, we will use some redundant formulation of the method for the transformed problem with three sets of coefficient matrices $A_n,B_n,K_n\in\R^{s\times s}$.
For the sake of an efficient implementation, the matrices $A_n$ and $K_n$ are lower triangular, and $A_n$ preferably has constant diagonal elements.
The additional index $n$ at e.g. $A_n=(a_{ij}^{(n)})_{i,j=1}^n$ indicates that in some very few time steps, especially at the boundaries, different coefficients may be used.
\par
In the following, by $e_n=(\delta_{nk})_k$ we denote the cardinal basis vectors in spaces of different dimensions, by $\eins=(1,\ldots,1)\T$ the vector of ones and by $I$ the identity matrix, the latter two sometimes with an additional index indicating the space dimension.
\par
Defining approximations $U_{ni}\approx u(t_n+c_ih_n),i=1,\ldots,s,$ and stacking the stage vectors into long vectors $Y_n=(Y_{ni})_{i=1}^s\in\R^{sm}$, $U_n=(U_{ni})_{i=1}^s\in\R^{sd}$, a two-step Peer method applied to the
ODE system \eqref{OCprob_ODE}-\eqref{OCprob_ODEinit} with constant step sizes is given by
\begin{align}
\label{PMv}
A_n Y_n=&\,B_nY_{n-1}+hK_nF(Y_n,U_n),\ n=1,\ldots,N.
\end{align}
with $F(Y_n,U_n):=\big(f(Y_{ni},U_{ni})\big)_{i=1}^s$.
This is already an abbreviated version since for the coefficient matrices like $A\in\R^{s\times s}$, we will use the same symbol for its Kronecker product $A\otimes I$ with the identity matrix as a mapping from the space $\R^{sm}$ to itself.
The abstract starting method $\Psi_s$ from \cite{SchroederLangWeiner2014} is specified now by an implicit Runge-Kutta method with one additional explicit term
\begin{align}
\label{IRKS}
A_0 Y_0=&\,a\otimes y_0+hb\otimes f(y_0,u_0)+hK_0F(Y_0,U_0),
\end{align}
with an appropriate approximation $u_0\approx u(0)$ and vectors
$a,b\in\R^s$.
Finally, the approximation for the solution $y(T)$ is slightly generalized by a linear combination
\begin{align}
\label{EndA}
y_h(T):=\sum_{i=1}^s w_i Y_{Ni}=(w\T\otimes I)Y_N,
\end{align}
where $\ones\T w=1$.
Of course, for ease of analysis and implementation, a standard method $(A,B,K)$ will be used for most of the time steps.
Besides the starting method \eqref{IRKS}, only the final step with $n=N$ will use a different method $(A_N,B_N,K_N)$ which needs to satisfy fewer order conditions only without harming the overall order.
\par
In the first-discretize-then-optimize approach the Lagrangian of the method has to be considered.
For the overall scheme \eqref{PMv}-\eqref{EndA} and multipliers $P=(P_0,\ldots,P_N)^T$, it is given by
\begin{align}
\notag
L(Y,P)
&=\,-C(y_h(T))\\
\notag
&\,+P_0\T\big(A_0 Y_0-a\otimes y_0-hb\otimes f(y_0,u_0)-hK_0 F(Y_0,U_0)\big)\\
\label{LsgFkG}
&\,+\sum_{n=1}^{N} P_n\T\big(A_nY_n-B_nY_{n-1}-hK_nF(Y_n,U_n)\big).
\end{align}
Computing the derivatives with respect to $Y_n$ leads to three different cases,
\begin{align}
\label{AdjEnd}
A_0\T P_0=&\,B_1\T P_1+
h\nabla_YF(Y_0,U_0)\T K_0\T P_0,\,n=0,\\
\label{AdjABK}
A_n\T P_n=&\,B_{n+1}\T P_{n+1}+h\nabla_YF(Y_n,U_n)\T K_n\T P_n,\ 0\le n\le N-1,\\
\label{AdjStrt}
A_N\T P_N=&\,w\otimes p_h(T)+h\nabla_YF(Y_N,U_N)\T K_N\T P_N,\ n=N.
\end{align}
Here, $p_h(T)=\nabla_y C\big(y_h(T)\big)\T$ and the Jacobian of $F$ is a block diagonal matrix $\nabla_YF(Y_n,U_n)=\diag_i\big(\nabla_{Y_{ni}}F(Y_{ni},U_{ni})\big)$.
Unfortunately, all these equations contain expressions of the form $\nabla_{Y_{j}}F\sum_{i}P_{i}k_{ij}$ which may be interpreted as a half-one-leg form for the adjoint right-hand side $\phi=-(\nabla_y f)\T p$.
Since such a scheme may be very difficult to analyze, we restrict the matrices $K_n,\,0\le n\le N,$ to diagonal form, which also means $K_n\T=K_n$.
Then, substituting the discrete controls $U_n=U_n(Y_n,P_n)$ in terms of
$(Y_n,P_n)$ and defining $\Phi(Y_n,P_n):=\big(\phi(Y_{ni},P_{ni})\big)_{i=1}^s$, the equations \eqref{AdjEnd}-\eqref{AdjStrt} can be rewritten as an approximation for the adjoint differential equation $p'=\phi(y,p)$ in the form
\begin{align}
\label{ADJEnd}
A_0\T P_0=&\,B_1\T P_1-hK_0\Phi(Y_0,P_0),\,n=0,\\
\label{ADJABK}
A_n\T P_n=&\,B_{n+1}\T P_{n+1}-hK_n\Phi(Y_n,P_n),\ 0\le n\le N-1,\\
\label{ADJStrt}
A_N\T P_N=&\,w\otimes p_h(T)-hK_N\Phi(Y_N,P_N),\ n=N.
\end{align}
With the restricted diagonal form of the matrices $K_n$, we still gain $s\!-\!1$ degrees of freedom per method compared to the simple case $K_n\!=\!\kappa I$, which has been considered in \cite{SchroederLangWeiner2014}.
\par
These equations are accompanied by the scheme \eqref{PMv}-\eqref{EndA}. Substituting once again the vectors $U_n$ by $(Y_n,P_n)$ results in the
following approximation for the forward equation $y'=g(y,p)$:
\begin{align}
\label{FORWABK}
A_n Y_n=&\,B_n Y_{n-1}+h K_n G(Y_n,P_n),\ n=1,\ldots,N,\\
\label{FORWStrt}
A_0 Y_0=&\,a\otimes y_0+hb\otimes g\big(y_0,p_h(0)\big)+hK_0G(Y_0,P_0),\\
\label{FORWEnd}
y_h(T)=&\,(w\T\otimes I)Y_N.
\end{align}
Similar to (\ref{EndA}), the value of $p_h(0)$ is determined by
an interpolant $p_h(0)=(v\T \otimes I)P_0$ of appropriate order
with $\eins\T v=1$.
\par
The key observation for the consistency analysis of the overall scheme
\eqref{ADJEnd}-\eqref{FORWEnd} is that these
discrete equations can be viewed as a discretization of the
two-point boundary value problem \eqref{RWPy}-\eqref{RWPp}.

\section{Order conditions}\label{SOrdng}
Order conditions for Peer methods are obtained by Taylor expansions of its residuals with the function values of the exact solutions $y$, resp. $p$.
Defining the partial sums $\exp_q(z):=\sum_{j=0}^{q-1} z^j/j!$ with $q$ terms, Taylor's theorem for a smooth function $v\in C^{q}[0,T]$ may be written as
\[ v(t_n+c_ih)=\exp_q(c_iz) v|_{t=t_n}+O(h^q),\ z:=h\frac{d}{dt},\]
with some slight abuse of notation.
Introducing the column vector $\cc=(c_i)_{i=1}^s\in\R^s$ of nodes, $\exp_q(\cc)\in\R^s$ is defined by component-wise application.
Expanding the residuals of \eqref{FORWABK} with the values of $y$ for order $q_1$ and \eqref{ADJABK} with values of $p$ for order $q_2$ gives
\begin{align}
\label{ordep}
\Big(A_n\exp_{q_1}(\cc z)-B_n\exp_{q_1}\big((\cc-\ones)z\big)-zK_n\exp_{{q_1}-1}(\cc z)\Big)y|_{t_n}\stackrel!=O(z^{q_1}),\\
\label{ordeq}
\Big(A_n\T\exp_{q_2}(\cc z)-B_{n+1}\T\exp_{q_2}\big((\cc+\ones)z\big)+zK_n\exp_{q_2-1}(\cc z)\Big)p|_{t_n}\stackrel!=O(z^{q_2}).
\end{align}
We note, that $(q_1,q_2)$ correspond to the local orders of the methods.
For a representation in matrix form, the Vandermonde matrix
\begin{align}\label{VdM}
 V_q(\cc):=\Big(\ones,\cc,\cc^2,\ldots,\cc^{q-1}\Big)\in\R^{s\times q}
\end{align}
and $V_q'(\cc):=\Big(0,\ones,2\cc,\ldots,(q-1)\cc^{q-2}\Big)\in\R^{s\times q}$ are introduced.
By the binomial formula, shifts of nodes correspond to multiplications of $V(\cc)$ by the upper triangular Pascal matrix $\PP_q=\big({j-1\choose i-1}\big)_{i,j=1}^q$ containing the binomial coefficients.
In fact, we have $V_q(\cc-\ones)=V_q(\cc)\PP_q^{-1}$ and $V_q(\cc+\ones)=V_q(\cc)\PP_q$.
Also, with the nilpotent matrix $\tilde E_q:=\big(i\delta_{i+1,j}\big)_{i,j=1}^q$, which commutes with the Pascal matrix since $\PP_q=\exp(\tilde E_q)$, it holds that $V_q'(\cc)=V_q(\cc)\tilde E_q$, see e.g. \cite{SchmittWeinerPodhaisky2005}.
Hence, the matrix versions of the order conditions \eqref{ordep}, \eqref{ordeq} are
\begin{align}
\notag
A_nV_{q_1}(\cc)=&\,B_n V_{q_1}(\cc-\ones)+K_nV_{q_1}'(\cc)\\
\label{ordvq1}
=&\,B_nV_{q_1}(\cc)\PP_{q_1}^{-1}+K_nV_{q_1}(\cc)\tilde E_{q_1},\\
\notag
A_n\T V_{q_2}(\cc)=&\,B_{n+1}\T V_{q_2}(\cc+\ones)-K_nV_{q_2}'(\cc)\\
\label{ordaq2}
=&\,B_{n+1}\T V_{q_2}(\cc)P_{q_2}-K_nV_{q_2}(\cc)\tilde E_{q_2}.
\end{align}
Comparing with \cite{SchroederLangWeiner2014} for nonsingular $K_n$, this means that the forward conditions of that paper apply to the method $(K_n^{-1}A_n,K_n^{-1}B_n,I)$ while the adjoint conditions apply to $(A_nK_n^{-1},B_{n+1}K_n^{-1},I)$.
So, indeed, the redundant formulation \eqref{PMv} introduces $s\!-\!1$ additional degrees of freedom.
\par
Since all versions of Vandermonde matrices have been reduced to $V_q(\cc)$, now we may drop the argument in the remaining text, $V_q:=V_q(\cc)$.
The forward and adjoint starting methods \eqref{FORWStrt} and \eqref{ADJStrt} are Runge-Kutta methods lacking the computation of a final solution.
Therefore,
it has to be ensured that the linear combinations
$p_h(0)=v\T\otimes P_0$ and $y_h(T)=w\T\otimes Y_N$ are
$O(h^s)$-approximations to $p(0)$ and $y(T)$, if $P_0$ and $Y_N$ are
$O(h^s)$-approximations to $p$ and $y$ themselves, respectively.
\begin{lemma}
\label{Lintp}
With some vectors $v,w\in\R^s$ and nodes $c_1,\ldots,c_s$, the identities
\begin{align*}
  \pi(0)=\sum_{i=1}^s v_i \pi(hc_i)\;\text{ and }\;
  \pi(h)=\sum_{i=1}^s w_i \pi(hc_i),\; h\in\R,
\end{align*}
hold for all polynomials $\pi$ of degree $s-1$ iff
\begin{align}
  \label{ordv}
  v\T \ones = 1,\;v\T\cc^j = 0,\ j=1,\ldots,s-1,\\
  \label{ordw}
  w\T \cc^j = 1,\ j=0,\ldots,s-1.
\end{align}
\end{lemma}
\noindent{\bf Proof:}
Straightforward by changing the order of summation.\qed
\par
In what follows, we will choose $v$ and $w$ accordingly.
The accuracy of the approximations $Y_0$ and $P_N$ are now determined
by the stage orders $(q_1,q_2)$ which are derived in a way analogous to \eqref{ordep}, \eqref{ordeq}.
Formally, the two order conditions for \eqref{FORWStrt} and \eqref{ADJStrt} are
\begin{align}
\label{ordsrt}
A_0\exp_{q_1}(\cc z)=&\,a+bz+zK_0\exp_{q_1-1}(\cc z)+O(h^{q_1}),\\\label{ordads}
A_N\T\exp_{q_2}(\cc z)=&\,w\exp(z)-zK_N\exp_{q_2-1}(\cc z)+O(h^{q_2}).
\end{align}
With $q_1\le s+1$ and $q_2\le s$, the matrix versions of these conditions follow as before:
\begin{align}
\label{ordsq1}
A_0V_{q_1}=&\,ae_1\T+be_2\T+K_0V_{q_1}\tilde E_{q_1},\\
\label{ordeq2}
A_N\T V_{q_2}=&\,w\ones\T-K_NV_{q_2}\tilde E_{q_2},
\end{align}
with the cardinal basis vectors $e_j\in\R^s,\,j=1,\ldots,s$.
For $s\ge q_1\ge2,\,s\ge q_2\ge1$ the properties $V_qe_1=\ones$, $\tilde E_qe_1=0$, $\tilde E_qe_2=e_1$ have the following simple consequences
\begin{align}\label{precons}
 a=A_0\ones,\ b=A_0\cc-K_0\ones,\ w=A_N\T\ones.
\end{align}
\par
Since the combined schemes require many different order conditions, for ease of reference the specific choices are listed in Table~\ref{TabOrd}.
\begin{table}
\begin{center}
\begin{tabular}{|l|c|c|}\hline
 Steps&forward&adjoint\\\hline
 Start, $n=0$&\eqref{ordsq1} with $q_1=s$&\eqref{ordaq2} with $q_2=s-1$\\
 Standard, $1\le n<N$&\eqref{ordvq1} with $q_1=s+1$&\eqref{ordaq2} with $q_2=s$\\
 Last step&\eqref{ordvq1}, $n=N$, $q_1=s$&\eqref{ordaq2}, $n=N-1$, $q_2=s-1$\\
 End point& \eqref{ordw}&\eqref{ordeq2} with $q_2=s-1$\\\hline
\end{tabular}
\caption{Combined order conditions for the different steps}\label{TabOrd}
\end{center}
\end{table}
\section{Convergence}\label{SKonvg}
In this section the errors $\check Y_{nj}:=y(t_{nj})-Y_{nj}$, $\check P_{nj}:=p(t_{nj})-P_{nj},$ $n=0,\ldots,N$, $j=1,\ldots,s$, are analyzed.
It is convenient to multiply the forward Peer steps by $A_n^{-1}$.
This gives new coefficient matrices $\bar B_n:=A_n^{-1}B_n$ and $\bar K_n:=A_n^{-1}K_n$.
For the general forward step \eqref{FORWABK}, we obtain the relation
\begin{align}\label{Fevorn}
 \check Y_n=\bar B_n\check Y_{n-1}
 +h\bar K_n(\nabla_y G_n \check Y_n+\nabla_p G_n \check P_n)+\tau_n^Y,
\end{align}
$1\le n\le N$.
Here, $\tau_n^Y$ ist the truncation error and the matrix derivatives are block diagonal matrices and placeholders for integral mean values as in
\[ g(Y_{nj}+\check Y_{nj},P_{nj})-g(Y_{nj},P_{nj})
 =\int_0^1 \nabla_y g(Y_{nj}+\sigma\check Y_{nj},P_{nj})d\sigma\cdot \check Y_{nj}
\]
for $\nabla_y G_n\cdot \check Y_n$.
In the starting step \eqref{FORWStrt}, $B_0$ is missing but there is an additional $O(h)$-contribution from $g(y_0,p_h(0))$,
\begin{align}\label{FeFStrt}
 \check Y_0=h\bar K_0\nabla_y G_0 \check Y_0
  +h(\bar K_0\nabla_pG_0+\bar bv\T\otimes \nabla_p g_0)\check P_0
  +\tau_0^Y,
\end{align}
where $\bar b=A_0^{-1}b$.
We remind that according to Table~\ref{TabOrd} the truncation errors $\tau^Y$ satisfy $\tau_n^Y=O(h^{s+1})$ only for $1\le n\le N-1$ while $\tau_0^Y,\tau_N^Y=O(h^s)$ has lower order.
\par
The adjoint step \eqref{ADJABK} is multiplied by $A_n^{-T}$ and with abbreviations $\tilde B_{n+1}\T:=(B_{n+1}A_n^{-1})\T$ and $\tilde K_n\T:=(K_nA_n^{-1})^T$, the error equation becomes
\begin{align}\label{Feadjn}
 \check P_n=\tilde B_{n+1}\T \check P_{n+1}-h\tilde K_n\T(\nabla_y \Phi_n\check Y_n+\nabla_p \Phi_n\check P_n)+\tau_n^P.
\end{align}
This equation holds for $0\le n\le N-1$, since \eqref{ADJEnd} corresponds to \eqref{ADJABK} with $n=0$.
In the adjoint starting step \eqref{ADJStrt} the boundary condition reads $p_h(T)=\nabla_y C(y_h(T))\T$ in detail with $y_h(T)=(w\T\otimes I)Y_N$.
Since $A_N^{-T}w=\eins$ by \eqref{precons} this adjoint step gives rise to the equation
\begin{align}\label{FeAStrt}
 \check P_N=\left((\eins w\T )\otimes \nabla_{yy}C_N\right)\check Y_N -h\tilde K_N\T(\nabla_y\Phi_N\check Y_N+\nabla_p\Phi_N\check P_N)+\tau_N^P,
\end{align}
again with a mean value $\nabla_{yy}C_N\in\R^{m\times m}$ of the symmetric Hessian matrix of $C$.
\par
By numbering the unknowns in the order $\check Y_0,\ldots,\check Y_N,\check P_0,\ldots,\check P_N$ and the equations likewise the error equations \eqref{Fevorn}--\eqref{FeAStrt} give rise to a linear system
\begin{align}\label{BlkMh}
 \MM_h\check Z=\tau,
 \ \check Z=\begin{pmatrix} \check Y\\\check P \end{pmatrix},
 \ \tau=\begin{pmatrix}\tau^Y\\\tau^P \end{pmatrix},
\end{align}
where $\MM_h$ has a $2\times2$-block structure.
The terms not depending on $h$ are critical with respect to stability.
Hence, we look closer at the matrix $\MM_0$ in which all $O(h)$-terms have been deleted.
The matrix $\MM_0$ has lower block triangular structure with Kronecker products
\begin{align}\label{BlockM}
 \MM_0=\begin{pmatrix}
 M_{11}\otimes I_m&0\\
 M_{21}\otimes\nabla_{yy}C_N&M_{22}\otimes I_m
\end{pmatrix},
\end{align}
where $M_{11},M_{21},M_{22}\in\R^{s(N+1)\times s(N+1)}$.
For convenience, the index range corresponds to that of the grid, the blocks of, e.g., the first matrix are $(M_{11})_{ij}\in\R^{s\times s},\,0\le i,j\le N$. Its inverse is given by
\begin{align}\label{BlockMi}
 \MM_0^{-1}=\begin{pmatrix}
 M_{11}^{-1}\otimes I_m&0\\
 -(M_{22}^{-1}M_{21}M_{11}^{-1})\otimes\nabla_{yy}C_N&M_{22}^{-1}\otimes I_m
\end{pmatrix},
\end{align}
It is obvious that the lower block $(\MM_0)_{21}\in\R^{s(N+1)\times s(N+1)}$ is trivial for linear objective functions $C$.
By \eqref{FeAStrt} it contains one nontrivial block only in the last diagonal block of size $(sm)\times(sm)$ and $M_{21}$ does only have rank one.
In fact
\begin{align}\label{Msubdi}
  (\MM_0)_{21} 
 = M_{21}\otimes\nabla_{yy}C_N,
 \quad
 M_{21}=(e_{N}\otimes\eins)(e_{N}\otimes w)\T,
\end{align}
where $e_N=\big(\delta_{Nk}\big)_{k=0}^N$.
The factors $M_{jj},\,j=1,2,$ of the diagonal blocks of $\MM_0$ have again block structure in bi-diagonal form with identity blocks $I_s$ in the diagonal.
However, $M_{11}$ has nontrivial block subdiagonals $(M_{11})_{n,n-1}=-\bar B_n,\,1\le n\le N$, while the second matrix has block superdiagonals $(M_{22})_{n,n+1}=-\tilde B_{n+1}\T,\,0\le n<N$.
As before, the coefficient matrices $A_n,B_n$, $1\le n\le N-1$, from the standard scheme do not depend on the index.
It is easy to compute the blocks $(M_{kk}^{-1})_{ij},\,k=1,2,$ of its inverses explicitly.
In order to prove the convergence result, we need some special norm bound for these inverses.
It is well known, e.g. \cite{HairerNorsettWanner2006}, that due to zero stability there exist nonsingular matrices $X_1,X_2\in\R^{s\times s}$ such that $\|A^{-1}B\|_{X_1}:=\|X_1^{-1}A^{-1}BX_1\|_\infty=1$ and $\|(BA^{-1})\T\|_{X_2}:=\|X_2^{-1}(BA^{-1})\T X_2\|_\infty=1$.
For convenience we assume that this holds for all nontrivial blocks, $\|\bar B_n\|_{X_1}=1$ and $\|\tilde B_n\T\|_{X_2}=1,\,1\le n\le N$.
This is no severe restriction, since all matrices have the right eigenvector $\eins$ and possible exceptions concern two indices at most and may only spoil the constants of the following results.
The vector norms belonging to these matrix norms are $\|Y_n\|_{X_1}=\|X_1^{-1}Y_n\|_\infty$ and $\|P_n\|_{X_2}=\|X_2^{-1}P_n\|_\infty$.
\par
One of the norm bounds for the inverses depends on the block sparseness of the pre-image.
To this end we use the following notations.
For grid vectors $W\T=(W_0\T,\ldots,W_N\T)\in\R^{(N+1)sk},\,k\in\{1,m\}$, the {\em block sparsity} is denoted by
\begin{align}\label{spars}
  (\# W):=\{\# n:\, W_n\not=0,0\le n\le N\}\in\N_0.
\end{align}
The norm definitions are also extended to grid vectors, e.g. $\|Y\|_{X_1}:=\max\{\|Y_n\|_{X_1}:\,0\le n\le N\}$.
\begin{lemma}\label{LMNorm}
For the inverses of $M_{jj},\,j=1,2,$ from \eqref{BlockM}, the following estimates with pre-image $W\in\R^{(N+1)s}$ hold,
\begin{align}\label{nrmspar}
 \|M_{jj}^{-1}W\|_{X_j}\le&\;(\#W)\cdot\|W\|_{X_j},\ j=1,2,\\\label{nrmspars}
 \|M_{22}^{-1}M_{21}M_{11}^{-1}W\|_{X_2}\le&\;\gamma(\#W)\cdot\|W\|_{X_1},
\end{align}
with constant $\gamma>0$ and for $j=1,2$, we have
\begin{align}\label{NmMW}
 \|M_{jj}^{-1}W\|_{X_j}\le&\;
 \max\{2\|W_0\|_{X_j},2\|W_N\|_{X_j},N\|W_n\|_{X_j},1\le n<N\},\\
 \|M_{22}^{-1}M_{21}M_{11}^{-1}W\|_{X_2}\le&\;
 \gamma\max\{2\|W_0\|_{X_1},2\|W_N\|_{X_1},N\|W_n\|_{X_1},1\le n<N\}.
 \end{align}
\end{lemma}
\noindent {\bf Proof:}
The inverses have identity matrices in its diagonal blocks and the remaining blocks are easily verified to be
\begin{align}\label{Invblk}
 (M_{11}^{-1})_{nk}=\bar B_n\cdots\bar B_{k+1},\,k<n,\quad
 (M_{22}^{-1})_{nk}=\tilde B_{n+1}\T\cdots \tilde B_{k}\T,\,k>n.
\end{align}
Hence, due to assumption we have $\|(M_{11}^{-1})_{nk}\|_{X_1}=1,\,k\le n$, and $\|(M_{22}^{-1})_{nk}\|_{X_2}=1,\,k\ge n$, which leads for the first block to
\begin{align}\label{MiW}
 \|M_{11}^{-1}W\|_{X_1}\le \max_{0\le n\le N}\sum_{k=0}^n\|W_k\|_{X_1}
  \le\big(\|W_0\|_{X_1}+\|W_N\|_{X_1}\big)+\big(\sum_{k=1}^{N-1}\|W_k\|_{X_1}\big).
\end{align}
Now, \eqref{nrmspar} is a trivial consequence and \eqref{NmMW} follows by treating the two brackets on the right separately.
For $\MM_{22}^{-1}$ analogous estimates hold.
\par
The rank-one structure from \eqref{Msubdi} leads to the representation
\begin{align*}
 M_{22}^{-1}M_{21}M_{11}^{-1}=\eins_{s(N+1)}\omega\T,
 \quad \omega\T=\eins_{N+1}\T\otimes(\eins\T A)
\end{align*}
which is seen in the following way.
By \eqref{Msubdi} in the column vector $M_{22}^{-1}(e_N\otimes \eins)$ only the last column of $M_{22}^{-1}$ contributes and by \eqref{Invblk} these contibutions are $\bar B_{n+1}\cdots\bar B_N\T\eins=\eins$ by \eqref{ordaq2}.
In a similar way in $(e_N\otimes w)\T M_{11}^{-1}$ the last row, see \eqref{Invblk}, contributes $w\T\bar B_N\cdots\bar B_{k+1}=w\T A_N^{-1} B_N\cdots\bar B_{k+1}=\eins\T B_N\cdots\bar B_{k+1}=\eins\T A$, again by \eqref{ordaq2}.
Hence, for the subdiagonal block the estimate \eqref{MiW} appears again with the additional factor $\gamma=\|X_2^{-1}\eins\eins\T A X_1\|_\infty$.
\qed
\par
\begin{remark}
Of course, these estimates carry over for the block matrix \eqref{BlockM} to higher dimensions $m>1$ with $W\in\R^{(N+1)sm}$, and norms $\|W\|:=\max\{\|W_n\|_{X\otimes I_m}:\,0\le n\le N\}$.
\end{remark}
\par
Writing $\MM_h=\MM_0-h\Uu$ in the error equation \eqref{BlkMh}, it may be rewritten in fixed-point form
\begin{align}\label{FeFP}
 \check Z=h\MM_0^{-1}\Uu\check Z+\MM_0^{-1}\tau.
\end{align}
Here, an important point is that the matrix $\Uu$ contains exactly two nontrivial blocks in each column, which becomes obvious after inspecting the error equations \eqref{Fevorn}--\eqref{FeAStrt}.
This means that  Lemma~\ref{LMNorm} may be applied to the first term on the right hand side of \eqref{FeFP} with $(\#\Uu\check Z)=2$ leading to an $O(h)$-contraction.
\begin{theorem}\label{TKonv}
Let the Peer method with $s>1$ stages satisfy the order conditions collected in Table~\ref{TabOrd} and let the solutions satisfy $y\in C^{s+1}[0,T]$, $p\in C^s[0,T]$.
Assume, that a Peer solution $(Y\T,P\T)\T$ exists and that $f$ and $C$ have bounded second derivatives.
Then, for stepsizes $h\le h_0$ the error of these solutions is bounded by
\begin{align}\label{Fbeide}
  \|Y_{nj}-y(t_{nj})\|_\infty,\|P_{nj}-p(t_{nj})\|_\infty=O(h^{s-1}),
 \end{align}
$n=0,\ldots,N,\,j=1,\ldots,s$.
\end{theorem}
\noindent {\bf Proof:}
As a first step we inspect the inhomogeneity in \eqref{FeFP}.
Due to the block structure \eqref{BlockM} of $\MM_0$, we have
\begin{align}\label{FeIH}
 \MM_0^{-1}\tau=\begin{pmatrix}
  (M_{11}^{-1}\otimes I_m)\tau^Y\\
  (\MM_0^{-1})_{21}\tau^Y+(M_{22}^{-1}\otimes I_m)\tau^P
 \end{pmatrix}
 =\begin{pmatrix} O(h^s)\\ O(h^{s-1}) \end{pmatrix}.
\end{align}
These orders are verified, e.g., for the first block with \eqref{NmMW} by the assumptions in Table~\ref{TabOrd} through
\[ \|(M_{11}^{-1}\otimes I_m)\tau^Y\|_{X_1}
 \le \gamma(2h^s+2h^s+Nh^{s+1})= O(h^s),
\]
with a generic constant $\gamma$.
Due to lower order requirements  the second block has order $h^{s-1}$ only.
\par
The second step has to show that \eqref{FeFP} is a contractive fixed-point equation.
Inspection of the error equations \eqref{Fevorn}--\eqref{FeAStrt} shows that in each block column the matrix $\Uu$ has exactly one block entry in the main diagonal and one block entry in the diagonal of the off-diagonal blocks of $\MM_h$.
Hence, we have
\begin{align}\label{DefWQ}
 \Uu\begin{pmatrix}\check Y\\\check P\end{pmatrix}
 =\begin{pmatrix} W\\Q \end{pmatrix}
\end{align}
with $(\# W)=(\# Q)=1$.
Due to the block triangular form of $\MM_0$ and by \eqref{nrmspar} it follows that
\begin{align*}
\|\MM_0^{-1}\Uu\begin{pmatrix}\check Y\\\check P\end{pmatrix}\|
 =&\max\{\|(M_{11}^{-1}\otimes I_m)W\|_{X_1},
  \|(\MM_0^{-1})_{21}W+(M_{22}^{-1}\otimes I_m)Q\|_{X_2}\}\\
  \le&\max\{\|W\|_{X_1},\gamma\|W\|_{X_1}+\|Q\|_{X_2}\}
  \le L\|\check Z\|.
\end{align*}
The constant $L$ contains bounds for the derivatives of $g$ and $\phi$.
Now, by the Banach fixed-point theorem, equation \eqref{FeFP} is uniquely solvable for $h\le h_0=1/(2L)$ and the solution is bounded by
\begin{align*}
 \|\check Z\|\le\frac1{1-hL}\|\MM_0^{-1}\tau\|
  \le& 2\max\{\|M_{11}^{-1}\tau^Y\|_{X_1},\|M_{22}^{-1}\tau^P\|_{X_2}\}
  = O(h^{s-1}),
\end{align*}
by \eqref{FeIH}. \qed
\par
\begin{remark}
The Theorem only applies to the methods \verb+Peer3o32w+ and \verb+BDF3o32+ from Sections \ref{SP3o32w} and \ref{SB3o32} below, since the semi-explicit end method in \verb+BDF3o22+ misses the order requirements.
The real matrix
\[ X_1=\begin{pmatrix}
 \frac13&\frac{41}{42}&\frac{-1}{12}\\[2mm]
 \frac13&\frac13&\frac{11}{42}\\[2mm]
 \frac13&\frac8{231}&\frac2{11}
\end{pmatrix}\]
transforms $\bar B$ from BDF3 \eqref{BDF3n} to real Jordan form and the norms for the end methods exceed 1 only slightly.
We have $\|\bar B_N\|_{X_1}\le1.02$ for \eqref{BDF3o32N} and $\|\bar B_N\|_{X_1}\le1.22$ for \eqref{ANimp}.
\end{remark}
\par
Of course, the error estimate \eqref{Fbeide} is not very satisfactory since it states $O(h^{s-1})$-convergence for $Y$ only.
But with this global estimate, the result may be improved by a better consideration of the lower triangular block structure in \eqref{FeFP}.
\begin{lemma}\label{LKonv}
Under the assumptions of Theorem~\ref{TKonv}, the error in the $Y$-variable is of order $s$, i.e.
 \begin{align}\label{FeY}
  Y_{nj}-y(t_{nj})=O(h^{s}),
  \ n=0,\ldots,N,\,j=1,\ldots,s.
 \end{align}
\end{lemma}
\noindent {\bf Proof:}
By Theorem~\ref{TKonv} the term $W$ in \eqref{DefWQ} also satisfies
\[ \|(M_{11}^{-1}\otimes I_m)W\|_{X_1}\le\|W\|_{X_1}=O(h^{s-1})\]
since $(\# W)=1$.
Considering now the $\check Y$ part of \eqref{FeFP} only it is seen that
\[ \check Y=h(M_{11}^{-1}\otimes I_m)W+(M_{11}^{-1}\otimes I_m)\tau^Y
 =h O(h^{s-1})+O(h^s)=O(h^s).
\]
by \eqref{FeIH}.
\qed
\par
\begin{remark}\label{RKonvp}
An analogous discussion for the $P$-errors may explain some observations in the numerical tests below.
Here, one gets
\[ \|\check P\|_\infty\le \gamma_Y h^s+\gamma_P h^{s-1}\]
which is of order $O(h^{s-1})$ only, of course.
However, if the constant $\gamma_Y$ is much larger than the truncation error $\gamma_P h^{s-1}$, the observed orders may range between $s\!-\!1$ and $s$.
\end{remark}

\section{Construction of 3-stage methods}\label{SKonstr}
In \cite{SchroederLangWeiner2014} the adjoint boundary condition \eqref{ordads} was identified as the essential bottleneck for higher order.
The reason becomes obvious after writing the step \eqref{AdjStrt} out for $s=3$ with triangular matrix $A_N$ and $K_N=\diag(\kappa_i^{(n)})$:
\begin{align}
 \label{enda3}
 a_{33}^{(N)}{P_{N3}}=&\;w_3 p_h(T)-h\kappa_3^{(N)}\phi(Y_{N3},P_{N3}),\\
 \label{enda2}
 a_{22}^{(N)}P_{N2}+a_{32}^{(N)}P_{N3}=&\;w_2 p_h(T)-h\kappa_2^{(N)}\phi(Y_{N2},P_{N2}),\\
 \label{enda1}
 a_{11}^{(N)}P_{N1}+a_{21}^{(N)}P_{N2}+a_{31}^{(N)}P_{N3}=&\;
 w_1 p_h(T)-h\kappa_1^{(N)}\phi(Y_{N1},P_{N1}).
\end{align}
Obviously, with $K_N=I$ and $A_N=A$ in \cite{SchroederLangWeiner2014}, the first equation was an $O(h)$-approximation only of the correct boundary condition $P_{N3}=p_h(T)$ if $c_3=1$.
But with different coefficients $A_N\not=A$ and the redundant formulation \eqref{PMv} there are now three detours around this obstacle:
\begin{enumerate}
\item
With the choice $\kappa_3^{(N)}=0$ and $a_{33}^{(N)}=w_3\not=0$ equation \eqref{enda3} gives the exact boundary condition for $c_3=1$ and \eqref{enda2} corresponds to an implicit Euler step with sufficiently high local order 2.
A consequence of this choice is that the method  $(A_N,B_N,K_N)$ contains an explicit end stage.
\item
If $c_3<1$ the scheme \eqref{enda3} for the solution $P_{N3}$ may be an implicit Euler step for $p(t_N+hc_3)$.
\item
The triangular form of $A$ may be dropped for $A_N\not=A$.
The overall computational effort for the solution of the boundary value problem \eqref{RWPy}, \eqref{RWPp} increases only marginally if only the end step(s) have higher computational effort.
Later on, a fast converging simplified Newton iteration with triangular $\tilde A_N$ for the method $(A_N,B_N,K_N)$ is derived.
\end{enumerate}
\subsection{Three-stage standard Peer method}\label{Sstdpm}
For the internal time steps with $1\le n\le N-1$, a fixed method $(A,B,K)$ will be used.
With triangular form of $A$ and diagonal form of $K$, this method has $6+9+3=18$ free parameters plus 3 nodes for $s=3$.
On the other hand, the order conditions \eqref{ordvq1} and \eqref{ordaq2} comprise $3(q_1+q_2)$ conditions.
So it seems that the sum of the local orders may be bounded by $q_1+q_2\le 6$ resp. 7.
However, these order conditions are not independent and solutions exist beyond this bound.
Some background information of these dependencies is collected in Section~\ref{Sfao}.
It will be seen that all order conditions from Table~\ref{TabOrd} can only be satisfied with lowered adjoint orders which still is sufficient for order $O(h^3)$-convergence in the $y$-variable by Lemma~\ref{LKonv}.
Therefore, also the standard method $(A,B,K)$ is discussed with the lowered local order requirements $(q_1,q_2)=(s+1,s)=(4,3)$, too.
\par
Accordingly, the forward condition \eqref{ordvq1} is applied with $q_1=s+1=4$.
It is known that the method $(A,B,K)$ is invariant under a common shift of the nodes $\cc$.
Hence, for the sake of a simpler representation the following differences are introduced
\begin{align}\label{KnDif}
 d_1:=c_2-c_1,\ d_3:=c_3-c_2,
\end{align}
which means that $c_1=c_2-d_1$, $c_3=c_2+d_3$, and ordered nodes $c_1<c_2<c_3$ correspond to positive differences $d_j>0$.
Since the order conditions apply simultaneously to $A,B$, and their transposes, it is difficult to derive closed-form algebraic solutions.
Instead the conditions have been solved by algebraic manipulation with Maple.
Doing so, it turned out that both conditions could be solved by explicit substitutions up to $q_1=4$ and $q_2=2$.
For $q_2=3$ this is still the case for one component, but the two remaining conditions consist of highly nonlinear rational expressions.
However, the nominator of both conditions is essentially the same polynomial of high degree.
This polynomial consists of the factor $d_1\kappa_{2}$ and
\begin{align}\label{pord43}
 Q(d_1,d_3):=&
 3(11d_1^2+18d_1d_3+7d_3^2)d_1d_3-15d_1^3-67d_1^2d_3\\\notag
 &-55d_1d_3^2-7d_3^3
 +5(3d_1^2+5d_1d_3 +d_3^2)+3(d_1+d_3)-3.
\end{align}
Only the cancellation of $Q(d_1,d_3)$ makes sense, the other factors lead to confluent nodes or trivial methods since $\kappa_2$ cancels as a common factor of all matrices. 
\par
The solution set
\begin{align}\label{Qdef}
 {\cal Q}:=\{(d_1,d_3):\,Q(d_1,d_3)=0\}
\end{align}
is non-empty and defines a curve in the $(d_1,d_3)$-plane consisting of several probably unconnected branches which will be discussed later on in Subsection~\ref{SScan}.
\par
An interesting subclass of methods is defined by $d_3=d_1$ with nodes in equal distances.
In this case,
\begin{align}\label{Qdd}
 Q(d_1,d_1)=(3d_1-1)(2d_1-1)(6d_1^2-3d_1-1)
\end{align}
is a polynomial of degree 4 having 4 real solutions.
These solutions are:
\begin{itemize}
\item
 $d_1=d_3=\frac13$ essentially yields the BDF3 method.
The BDF3 method even fulfills the adjoint conditions \eqref{ordaq2} up to local order $q_2=4$, see \cite{SchroederLangWeiner2014}, and is $A(\alpha)$-stable with $\alpha=86.032^o$.
\item
for $d_1=d_3=\frac12$ local order $q_1=4$ is possible with $\kappa_1=0$ only leading to a blind first stage $Y_{n,1}=Y_{n-1,3}$.
After its elimination the BDF3 method is obtained again with larger stepsize.
\item
$d_1=d_3=\frac14-\frac{\sqrt{33}}{12}\approx -0.2287$. This method is not zero-stable, $\varrho(M(0))>1$ 
for the stability matrix \eqref{StbMat}.
\item
$d_1=d_3=\frac14+\frac{\sqrt{33}}{12}\approx 0.7287$. The method is $A(\alpha)$-stable with slightly larger $\alpha=87.871^o$ compared to BDF3.
However, since $c_3-c_1=d_1+d_3>1$, implementation of this method may be slightly more complicated.
\end{itemize}
\subsection{Scanning the parameter set ${\cal Q}$}\label{SScan}
For the Dahlquist test equation $y'=\lambda y$, one step of a Peer method $(A,B,K)$ reduces to a simple multiplication of the stage vector $Y_n$ by the stability matrix
\begin{align}\label{StbMat}
 M(z)=(A-zK)^{-1}B,\ z=h\lambda\in\C.
\end{align}
Zero-stability is the minimal requirement for a practical method and it means that the sequence $\big(M(0)^n\big)_{n\ge0}$ is bounded and requires that the eigenvalues of $M(0)$ lie inside the unit disc and those on the unit circle are semi-simple.
\par
\begin{figure}[t]
\centering
\begin{tabular}{cc}
\includegraphics[width=5.6cm]{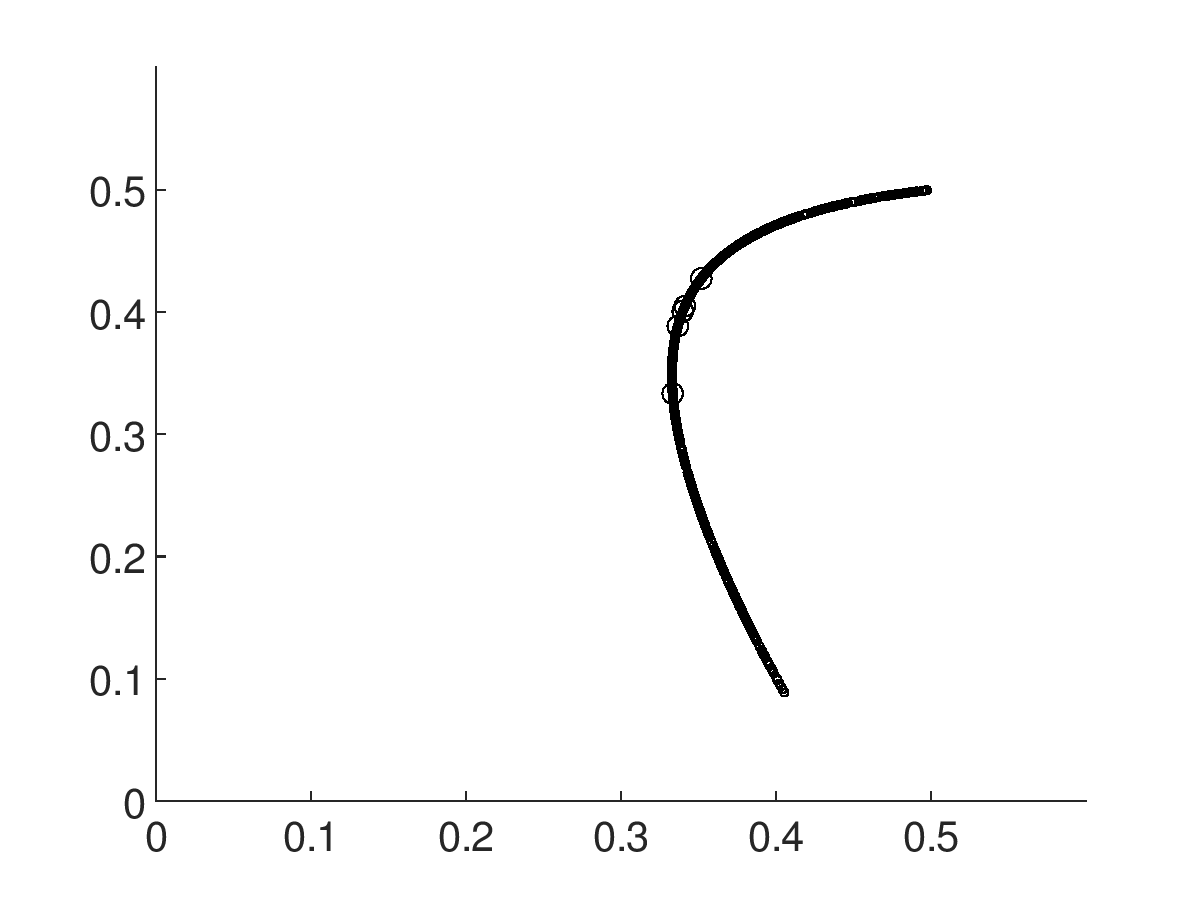}&
\includegraphics[width=5.6cm]{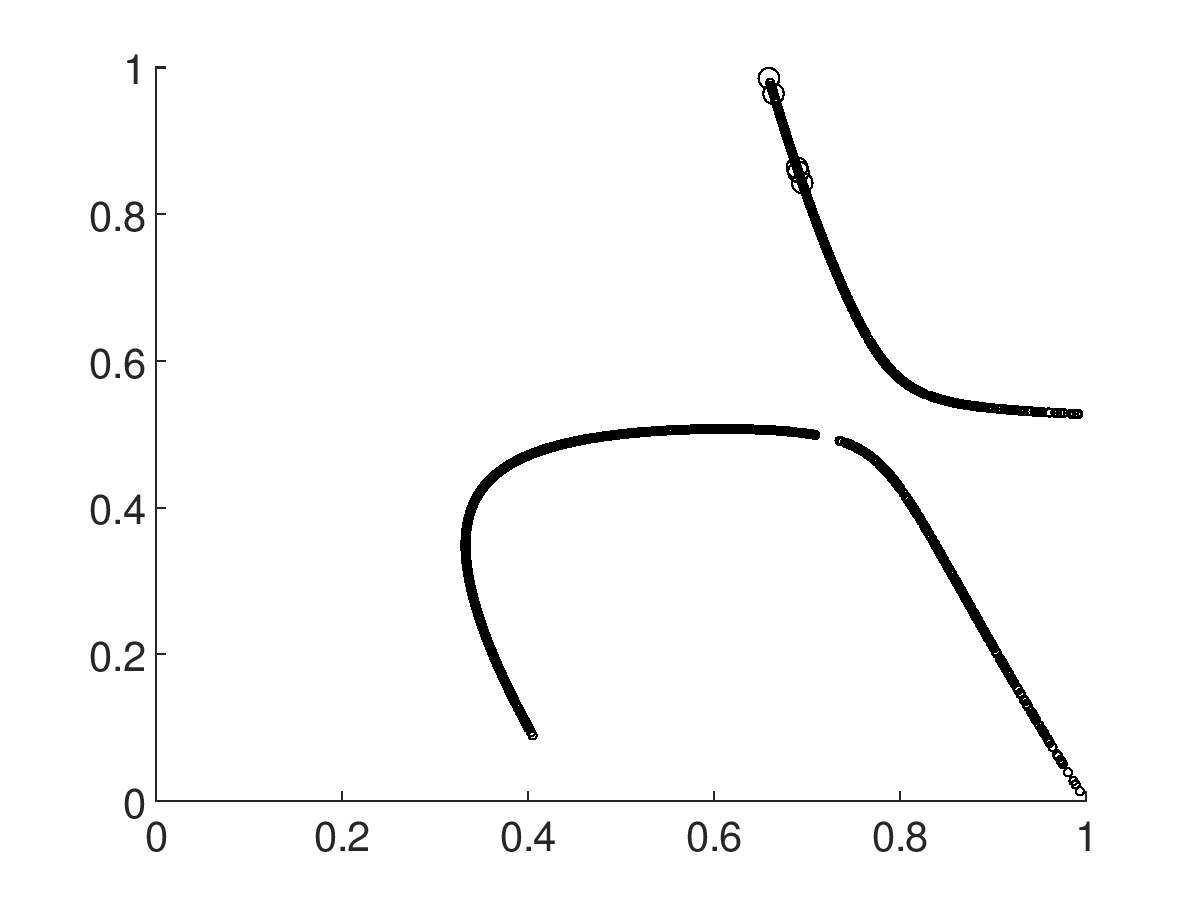}\\
\includegraphics[width=5.6cm]{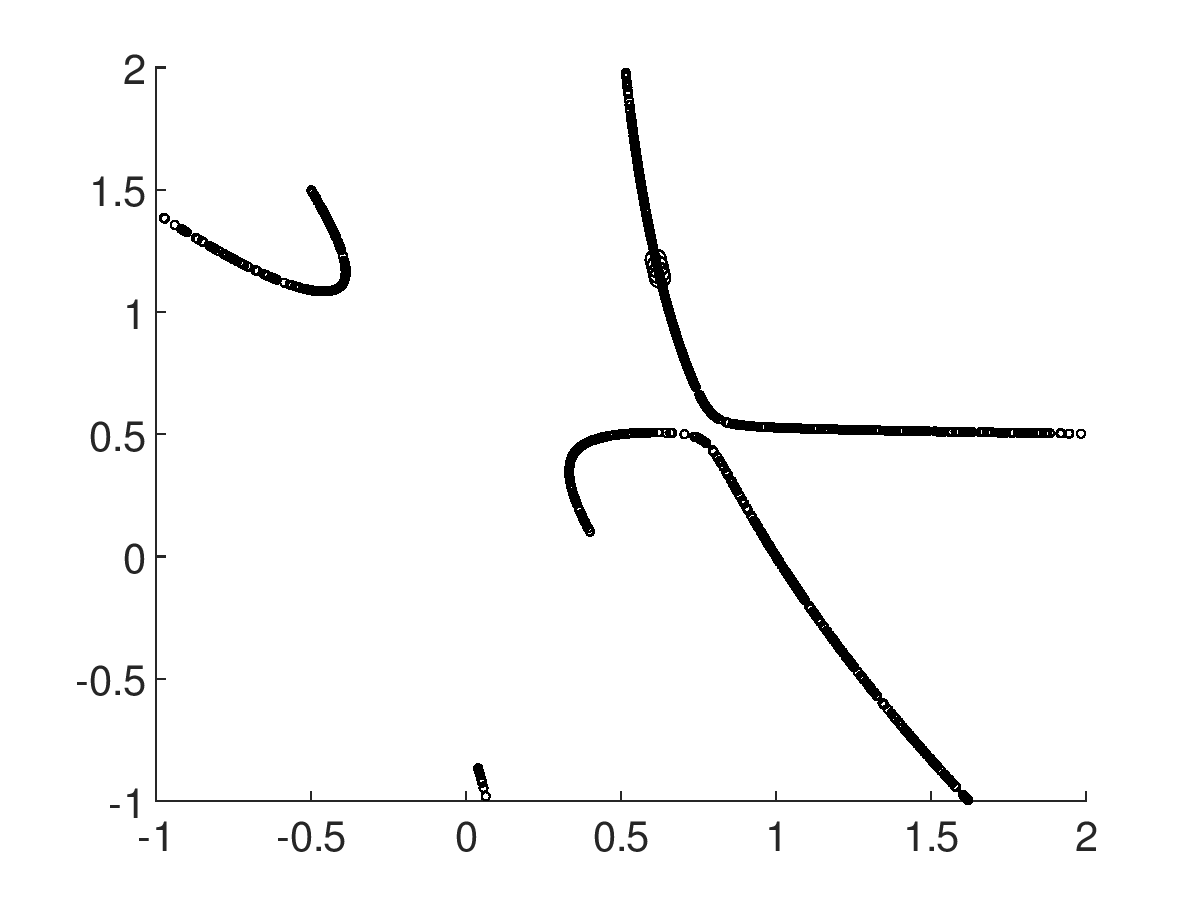}&
\includegraphics[width=5.6cm]{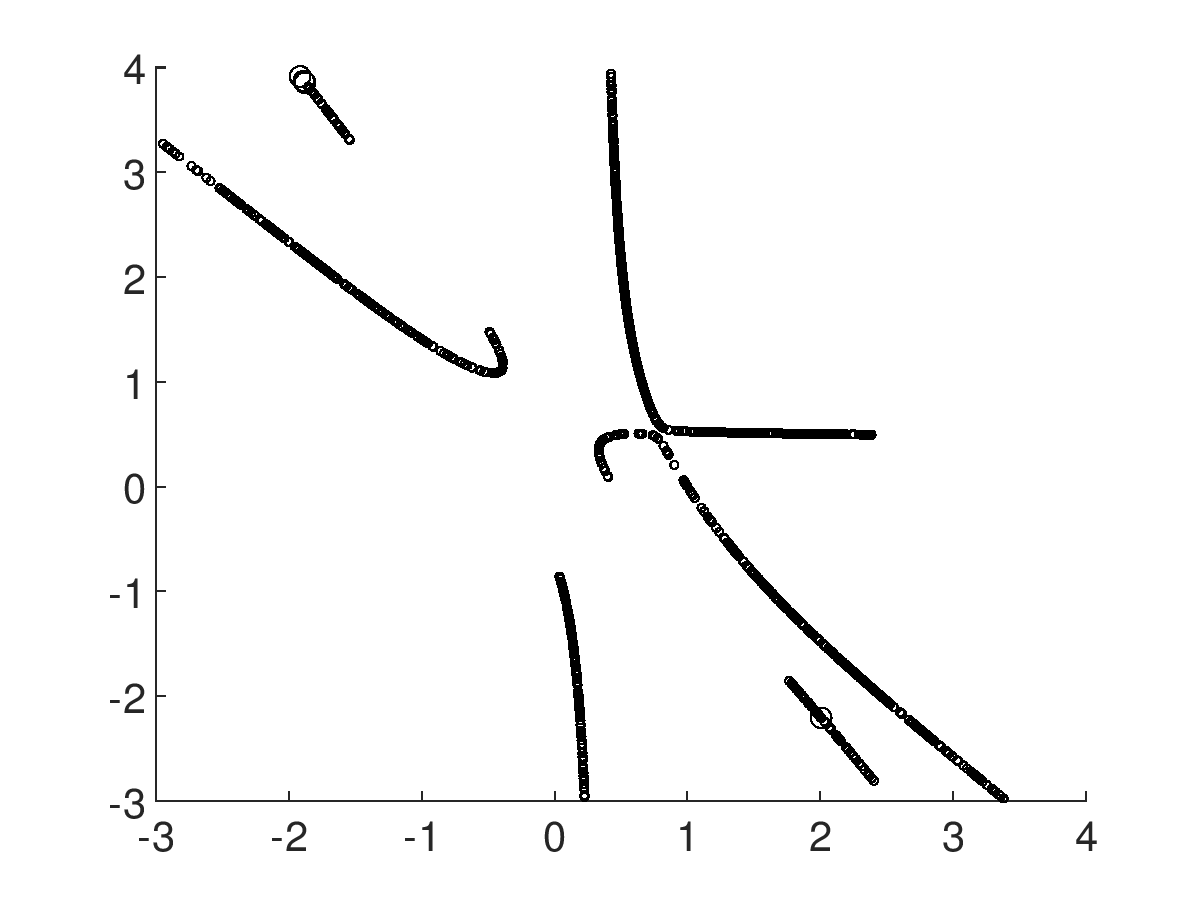}
\end{tabular}\\
\parbox{13cm}{
\caption{Sets in the $(d_1,d_3)$-plane with zero-stable methods for different zooms, nearly maximal stability angles marked with larger circles, see Table~\ref{Tabalf}.}
\label{Figalf}
}
\end{figure}

Of practical interest for stiff equations is $A(\alpha)$-stability which essentially means that
\begin{align}
\varrho\big(M(z)\big)<1\mbox{ for }\arg(z)<180^o-\alpha,
\end{align}
for the spectral radius $\varrho$.
Hence, $\big(M(z)^n\big)_n$ is bounded in a sector centered at the negative real axis with aperture $2\alpha$ and $A$-stability corresponds to $A(90^o)$-stability.
As for multistep methods, the corresponding angle $\alpha$ may be computed quite simply by reformulating the eigenvalue problem for $M(z)$ with some vector $x\in\C^s$:
\begin{align}\label{EWWOK}
 M(z)x=\lambda x\gdw Bx=\lambda(A-z K)x\gdw
 K^{-1}(A-\lambda^{-1}B)x=zx.
\end{align}
Solving the last equation as an eigenvalue problem for $z\in\C$ with $|\lambda|=|\lambda^{-1}|=1$ on the unit circle gives the root-locus-curves defining the boundary of the $A(\alpha)$-stability set.
Nearly maximal angles for some sets in the $(d_1,d_3)$-plane were computed in Matlab by starting a Gauss-Newton method for 2000 random points and computing some point from ${\cal Q}$ nearby.
After checking zero stability there, the maximal argument of eigenvalues $z$ of \eqref{EWWOK} were computed with 2000 points $\lambda$ on the unit circle.
The maximal angles and corresponding parameters are shown in Table~\ref{Tabalf}. The diagrams in Figure~\ref{Figalf} sketch those parts of the set ${\cal Q}$ belonging to zero-stable Peer methods for different zooms.
Larger circles in these diagram mark points with (nearly) maximal angles.
\begin{table}
\centering
\begin{tabular}{|l|c|c|c|c|}\hline\rule{0mm}{5mm}\hspace{-0.1cm}
 Set $(d_1,d_3)\in$& ${0\le d_1,d_3,\atop d_1+d_3\le1}$&$[0,1]^2$&$[-1,2]^2$&$[-3,4]^2$
 \\[2mm]\hline\rule{0mm}{4mm}\hspace{-0.1cm}
 Angle $\alpha$&86.194&88.341&88.419&90\\
 At $(d_1,d_3)$&$(0.3397,0.4)$&$(0.657,0.996)$&$(0.623,1.16)$&\\
 Diagram&top-left&top-right&lower left&lower right\\ \hline
\end{tabular}\\
\parbox{13cm}{
\caption{Maximal $A(\alpha)$-angles for some sets in the $(d_1,d_2$)-plane.}
\label{Tabalf}
}
\end{table}
These data show that for nodes $c_i$ in the standard interval $[0,1]$, the angle of BDF3 can be improved only marginally.
Nodes outside the interval $[0,1]$ may be less convenient but do not lead to difficulties and have been used before with Peer methods, see \cite{SchmittWeinerPodhaisky2005}.
Hence, the second method in Table~\ref{Tabalf} with a spread of $c_3-c_1\approx 1.65$ may be attractive in some cases since it more than halves the gap to A-stability.
The improvement for nodes in $[-1,2]^2$ is again marginal.
However, in $[-3,4]^2$ even A-stable methods seem to exist far out with rather exotic nodes, e.g. $(d_1,d_3)\approx (2.31,-2.66)$ or $(d_1,d_3)\approx (-1.56,3.34)$.
\par
Since, regarding stability, only a slight improvement over BDF3 is possible with nodes in $[0,1]$, a different criterion of practical interest is the leading error term.
For the method $(A,B,K)$ the forward error is given by
\begin{align}\label{Feq}
 \eta_{q+1}:=A\cc^{q+1}-B(\cc-\ones)^{q+1}-(q+1)K\cc^q
\end{align}
with $q=q_1$.
However, a similar scan as for the $A(\alpha)$-stability reveals that BDF3 has the minimal norm $\|\eta_5\|_\infty$ of all methods on ${\cal Q}$.
And since BDF coincides with its adjoint method and, hence, satisfies the order conditions with $q_2=q_1=s+1$ it is the first candidate for the standard method.
\subsection{Three-stage end methods with $\kappa_s^{(N)}=0$}\label{SKsing}
Setting $c_3=1$ and $w=e_3$, the order conditions \eqref{ordvq1}, \eqref{ordaq2} and \eqref{ordeq2} for $(A_N,B_N,K_N)$ can be easily fulfilled up to $q_1\!=\!q_2\!=\!2$.
But the additional condition \eqref{ordvq1} for the forward scheme with $q_1\!=\!3$ leads to a similar situation as in Section~\ref{Sstdpm}: the final conditions are multiples of a different polynomial $Q_N(d_1,d_3)$.
The combined conditions $Q(d_1,d_3)=Q_N(d_1,d_3)=0$ have several solutions which have to be computed numerically.
Unfortunately, all but one solution seem to lead to unfeasible methods with negative entries in $K$.
And the only feasible solution with $(d_1,d_3) \approx (0.417,0.0628)$ is not zero stable.
Hence, we have to work with an end method of local order $(2,2)$ only.
An example is given for the rescaled BDF3 method, $\cc\T=(\frac13,\frac23,1)$,
\begin{align}\label{BDF3n}
 A=\begin{pmatrix}
 \frac{11}6&0&0\\[1mm] -3&\frac{11}6&0 \\[1mm]\frac32&-3&\frac{11}6
 \end{pmatrix},
 \ B=\begin{pmatrix}
 \frac13&-\frac32&3\\[1mm] 0&\frac13&-\frac32\\[1mm] 0&0&\frac13
 \end{pmatrix},
 \ K=\frac13I.
\end{align}
The final method is
\begin{align}
 A_N=\begin{pmatrix}
 \frac{21}8&0&0\\[1mm] -\frac{14}3&\frac{23}{12}&0\\[1mm]
 \frac{49}{24}&-\frac{23}{12}&1
 \end{pmatrix},
 \ B_N=\begin{pmatrix}
 \frac12&-\frac{73}{24}&\frac{31}6\\[1mm] -\frac13&\frac{41}{12}&-\frac{35}6\\[1mm] \frac16&-\frac{37}{24}&\frac52\end{pmatrix},
\ K_N=\begin{pmatrix}\frac7{36}\\&\frac{23}{36}\\&&0
\end{pmatrix}.
\end{align}
An appropriate starting method with local orders $q_1=3$, $q_2=2$ is given with the coefficient matrices
\begin{align}\label{BDF3s}
 A_0=\begin{pmatrix}
 2& 0& 0\\[1mm] -\frac{10}3 &\frac{15}8&  0\\[1mm]
 \frac53& -\frac{73}{24}& \frac{11}6
 \end{pmatrix},
 \ K_0=\begin{pmatrix}
 \frac13&&\\[1mm] & \frac{25}{72}&\\[1mm] && \frac13
\end{pmatrix}.
\end{align}
This method will be called \verb+BDF3o22+ in the numerical tests.
\subsection{Three-stage end methods with $c_s\not=1$}\label{SP3o32w}
As mentioned before, for $c_s<1$ the adjoint boundary condition \eqref{enda3} is an implicit Euler step for $p(t_N+c_sh)$ and may be accurate enough for $s=3$, at least.
Now, condition \eqref{ordw} from Lemma~\ref{Lintp} becomes important.
As before, the conditions \eqref{ordvq1}, \eqref{ordaq2} and \eqref{ordeq2} for $(A_N,B_N,K_N)$ are easily solved for $q_1=q_2=2$.
For $q_1=3$ again only one condition remains which cannot be solved explicitly, a polynomial condition $Q_N(c_2,d_1,d_2)=0$.
Since it depends on the additional parameter $c_2$, any solution $(d_1,d_3)\in{\cal Q}$ may be plugged in and the equation is solved for $c_2$.
\par
For equidistant nodes $d_3=d_1$, the most interesting case in ${\cal Q}$ is $d_1=d_3=\frac13$ related to BDF3.
Here, the polynomial
\[ Q_N\left( c_2,\frac13,\frac13 \right)=12c_2^3-33c_2^2+28c_2-\frac{43}{6}\]
possesses three real roots $c_2\approx 0.48,\,0.92,\,1.3$.
Since the smallest value leads to nodes in the standard interval $[0,1]$, the corresponding end and starting methods are displayed.
The coefficients \eqref{BDF3n} of the standard method are invariant under node shifts and remain.
The entries in the boundary methods are polynomials in $c_2$ of degree 4 or less.
Hence, it is more convenient to give numerical values with
$c_2=0.48059993107999468110$, $c_1=c_2-1/3$ and $c_3=c_2+1/3$.
Denoting $A_N=(a_{ij}^{(N)})$, $B_N=(b_{ij}^{(N)})$, and $K_N=(k_{ij}^{(N)})$,
the coefficients of the final method are
\begin{equation}\label{BDF3o32N}
\begin{array}{lrlr}
 a_{11}^{(N)}=& 2,& b_{11}^{(N)}=&  0.5271726507800490190,\\[1mm]
 a_{12}^{(N)}=& 0,& b_{12}^{(N)}=& -2.0724604020801301580,\\[1mm]
 a_{13}^{(N)}=& 0,& b_{13}^{(N)}=&  3.5452877513000811390,\\[1mm]
 a_{21}^{(N)}=& -3.2608729312532042110,& b_{21}^{(N)}=& -0.3876786348934308516,\\[1mm]
 a_{22}^{(N)}=&  1.7608729312532043906,& b_{22}^{(N)}=&  1.4782541374935927700,\\[1mm]
 a_{23}^{(N)}=& 0,& b_{23}^{(N)}=& -2.5905755026001617388,\\[1mm]
 a_{31}^{(N)}=&  1.6957667700466743694,& b_{31}^{(N)}=&  0.19383931744671510930,\\[1mm]
 a_{32}^{(N)}=& -3.1888608156001606791,& b_{32}^{(N)}=& -0.57246040208012921227,\\[1mm]
 a_{33}^{(N)}=&  1.9930940455534862169,& b_{33}^{(N)}=&  0.87862108463341401017,\\
\end{array}
\end{equation}
and
\begin{align}
\notag
k_{11}^{(N)}=&\;0.32729496649332262670,\\[1mm]
\notag
k_{22}^{(N)}=&\;0.32125659965331187900,\\[1mm]
k_{33}^{(N)}=&\;0.37084850277337088940.
\end{align}
The starting method is
\begin{align}
 A_0=\begin{pmatrix}
  2.1796087544459576670&0&0\\[1mm]
 -4.2110754936961070457& 1.9644965156719027025&0\\[1mm]
 2.3648000725834827177& -3.1311631823385693702&\frac{11}6
 \end{pmatrix}
\end{align}
and
\begin{align}\label{BDF3o32s}
K_0=\begin{pmatrix}
 0.16049178284304720811&&\\[1mm]
 &0.37705439411285645618&\\[1mm]
 &&\frac13
 \end{pmatrix}.
\end{align}
This method is denoted by \verb+PEER3o32w+.
\subsection{Three-stage end methods with full $A_N$}\label{SB3o32}
In the setting of Section~\ref{SKsing} order $q_1=3$ of the end method could not be achieved.
The situation changes if we sacrifice the triangular form of $A_N$.
Of course, this increases the computational cost but since it concerns only one single time step the increase may be small compared to the overall cost.
In fact, an efficient iteration scheme based on a triangular matrix $\tilde A_N$ will be provided.
The order conditions  \eqref{ordvq1}, \eqref{ordaq2} and \eqref{ordeq2} for $(A_N,B_N,K_N)$ have no solution for $q_1=q_2=3$.
Solutions only exist for $q_1=3$, $q_2=2$ with $a_{1j}^{(N)},\,j=1,2,3,$ as free parameters.
Choosing block structure with $a_{12}^{(N)}=a_{13}^{(N)}=0$ in order to facilitate the construction of an iteration method and $a_{11}^{(N)}=\frac95$ for zero stability of the end step, the following end scheme for BDF3 is obtained:
\begin{align}\label{ANimp}\hspace{-0.35cm}
 A_N=\begin{pmatrix}
 \frac95& 0& 0\\[1mm]
 -\frac{109}{40}& \frac43& \frac7{24}\\[1mm]
 \frac{37}{40}& -\frac43&\frac{17}{24} \end{pmatrix},
 \ B_N=\begin{pmatrix}
 \frac{39}{80}& -\frac{19}{10}& \frac{257}{80}\\[1mm]
 -\frac{37}{120}&\frac{17}{15}& -\frac{77}{40}\\[1mm]
 \frac{37}{240}& -\frac25 & \frac{131}{240}
 \end{pmatrix},
 \ K_N=\begin{pmatrix}
 \frac7{24}&&\\[1mm]
 &\frac49&\\[1mm]
 &&\frac7{72}
 \end{pmatrix}.
\end{align}
An appropriate starting method uses
\begin{align}\label{A0imp}
 A_0=\begin{pmatrix}
 2& 0&0\\[1mm]
 -\frac{10}3&\frac{15}8&0\\[1mm]
 \frac53&-\frac{73}{24}&\frac{11}6
\end{pmatrix},
 \ K_0=\begin{pmatrix}
 \frac13&&\\[1mm]
 &\frac{25}{72}&\\[1mm]
 &&\frac13
 \end{pmatrix}.
\end{align}
The name of this method will be \verb+BDF3o32+.
\par
A simple implementation for the last time step \eqref{PMv}, $n=N$, is possible with a simplified Newton method where $A_N$ in the Jacobian is replaced by a lower triangular approximation $\tilde A_N$.
For the test equation $y'=\lambda y$ such an iteration has the form
\begin{align}\label{DrIter}
 (\tilde A_N-zK_N)(Y_N^{[1]}-Y_N^{[0]})=-A_NY_N^{[0]}+zK_NY_N^{[0]}+B_NY_{N-1},
\end{align}
$z=h\lambda$, and may be solved stage-by-stage.
The iteration matrix is $S(z):=(\tilde A_N-zK_N)^{-1}(\tilde A_N-A_N)$.
Moreover, with equal subdiagonals $\tilde a_{ij}^{(N)}=a_{ij}^{(N)},\,j<i,$ the stages $Y_{N1}^{[0]},Y_{N2}^{[0]},\ldots,$ may be overwritten and only slight modifications of the triangular forward step are necessary since $\tilde A_N-A_N$ has vanishing subdiagonals.
A good choice is
\begin{align}
\tilde A_N=\begin{pmatrix}
 \frac95& 0& 0\\[1mm]
 -\frac{109}{40}& \frac{73}{39}& 0\\[1mm]
 \frac{37}{40}& -\frac43&\frac{535}{752}
 \end{pmatrix},
\end{align}
which gives a very good contraction $\rho(S(z))\le 0.05$ for $z$ on the
negative real axis.
\par
In the adjoint boundary condition \eqref{ADJStrt}, an analogous procedure may be used and the convergence analysis applies as well.
The iteration matrix there is $(\tilde A_N\T+\zeta K_N)^{-1}(\tilde A_N-A_N)\T$ and it has the same eigenvalues as $S(-\zeta)$.

\section{Combined order conditions and symmetric nodes}\label{Sfao}
In \cite{SchroederLangWeiner2014} it was observed that the order conditions may simplify for nodes which are symmetric to some center point $\zeta/2\in\R$, which means that $\Pi\cc=\zeta\eins-\cc$ with the flip permutation $\Pi=\big(\delta_{i,s+1-j}\big)_{i,j=1}^s$ which is an involution, $\Pi^2=I$, and symmetric, $\Pi\T=\Pi$.
In this case the adjoint order conditions \eqref{ordaq2} correspond to the forward order conditions \eqref{ordvq1} with $q_2=q_1$ for the permuted matrices $\Pi K^{-1}A\T \Pi$ and $\Pi K^{-1} B\T \Pi$.
And since the original coefficients for the BDF method have Toeplitz form and, hence, are persymmetric, $\Pi K^{-1} A=A\T K^{-1}\Pi$, $\Pi K^{-1} B=B\T K^{-1}\Pi$, the adjoint conditions are satisfied automatically.
This question will be discussed now in more detail.
\par
From the previous discussions it seems that methods with nodes symmetric to $c_2$, i.e. $d_1=d_3$, may have superior convergence properties.
In this section we will look for possible reasons for that especially for higher orders $q_1,q_2\ge s$.
Rewriting the two order conditions \eqref{ordvq1} and \eqref{ordaq2} slightly as
\begin{align*}
AV_{q_1}\PP_{q_1}=&\,BV_{q_1}+KV_{q_1}\tilde E_{q_1}\PP_{q_1},\\
V_{q_2}\T A =&\,P_{q_2}\T  V_{q_2}\T B-\tilde E_{q_2}\T V_{q_2}\T K,
\end{align*}
and subtracting, after multiplying the first by $P_{q_2}\T  V_{q_2}\T$ from the left and the second by $V_{q_1}$ from the right, cancels $B$ and leaves the equation
\begin{align*}
 0=&\,P_{q_2}\T  V_{q_2}\T AV_{q_1}\PP_{q_1}-P_{q_2}\T  V_{q_2}\T KV_{q_1}\tilde E_{q_1}\PP_{q_1}-V_{q_2}\T AV_{q_1} -\tilde E_{q_2}\T V_{q_2}\T KV_{q_1}.
\end{align*}
This leads to the following Lemma.
\begin{lemma}\label{LSylv}
For any Peer method $(A,B,K)$ satisfying the order conditions \eqref{ordvq1} and the adjoint order conditions \eqref{ordaq2} with $q_1,q_2\in\N$ the matrices $A$ and $K$ are related by the following Sylvester-type matrix equation
\begin{align}\label{Sylvg}
(V_{q_2}\PP_{q_2})\T  A(V_{q_1}\PP_{q_1})-V_{q_2}\T AV_{q_1}
=&(V_{q_2}\PP_{q_2})\T KV_{q_1}\PP_{q_1}\tilde E_{q_1}+(V_{q_2}\tilde E_{q_2})\T KV_{q_1}.
\end{align}
\end{lemma}
We note that the operator on the left acting on $A$ is singular since $P_{q_1},P_{q_2}$ have 1 as a multiple eigenvalue.
More structure can be seen in equation \eqref{Sylvg} in the case $q_1=q_2=s$ where $V_s$ is non-singular.
After the congruence multiplication $(\cdots)\to V_s^{-T}(\cdots)V_s^{-1}$ two well-known matrices appear as coefficients.
The first one is the extrapolation matrix $\Theta:=V_sP_sV_s^{-1}$ and  $E:=V_s\tilde E_sV_s^{-1}$ is the differentiation matrix with respect to the nodes $\{c_i\}_{i=1}^s$.
Then $V_s\tilde E_sP_sV_s^{-1}=E\Theta=\Theta E$ and \eqref{Sylvg} is equivalent with
\begin{align}\label{Sylvs}
 \Theta\T A\Theta-A= \Theta\T K\Theta E+E\T K.
\end{align}
\subsection{Symmetric nodes}
Symmetric nodes lead to special properties of the extrapolation and differentiation matrices $\Theta$ an $E$ in \eqref{Sylvs}.
\begin{lemma}\label{LSymKn}
Let the nodes be symmetric to some center point, $\Pi\cc=\zeta\ones-\cc$, $\zeta\in\R$, with the flip permutation satisfying $\Pi=\Pi\T$, $\Pi^2=I$.
Then, the following identities hold:
\begin{align}\label{PiVq}
 \Pi V_q=&V_q\Delta_qP_q^\zeta,\quad \Delta_q:=\diag(1,-1,\pm 1,\ldots)\in\R^{q\times q},\quad q\in\N,\\\label{PiThPi}
 \Pi\Theta\Pi=&\Theta^{-1},\quad \Pi E\Pi=-E,\quad q=s.
\end{align}
\end{lemma}
\noindent{\bf Proof:}
By the binomial formula after a shift of the nodes $\cc\mapsto\cc+\zeta\eins$, the Vandermonde matrix $V_q$ is multiplied from the right by $P_q^\zeta=\exp(\zeta\tilde E_q)$.
The trivial identity $\Delta_q\tilde E_q=-\tilde E_q\Delta_q$ shows that $\Delta_q\PP_q=\PP_q^{-1}\Delta_q$.
Considering $\Pi V_q$ column-wise for $1\le j\le q$, \eqref{PiVq} is the matrix version of the identity
\begin{align*}
 \Pi\cc^{j-1}=(\zeta\eins-\cc)^{j-1}
 =\sum_{i=1}^{j-1}\cc^{i-1}(-1)^{i-1}\underbrace{{j-1\choose i-1}\zeta^{j-i}}_{\PP^\zeta}.
\end{align*}
And this immediately yields
\begin{align*}
 \Pi\Theta\Pi = \Pi V_sP_s(\Pi V_s)^{-1}
  =V \Delta_s P_s^\zeta P_s P_s^{-\zeta}\Delta_s V^{-1}
  =V\Delta_s P_s\Delta_s V^{-1}=\Theta^{-1},
\end{align*}
since the Pascal matrix with checkerboard sign changes is its inverse.\qed
\par
A direct consequence is
\begin{lemma}
If there exist pairs $(A,K)$ solving the Sylvester equation \eqref{Sylvg} for $q_1=q_2$ and the nodes are symmetric, i.e. $\Pi\cc=\zeta\ones-\cc$, $\zeta\in\R$, then there is also a persymmetric solution pair with $\Pi A\T\Pi=A$, $\Pi K\Pi=K$.
\end{lemma}
\noindent {\bf Proof:}
Since both orders are equal, the index on $q_1=q_2$ may be dropped.
With the permutation $\Pi$, $\Pi^2=I$, by Lemma~\ref{LSymKn} the left hand side of \eqref{Sylvg} may be rewritten as
\begin{align*}
 &(V_{q}P_{q})\T  A(V_{q}\PP_{q})-V_{q}\T AV_{q}=
 \,(\Pi V_{q}P_{q})\T \Pi A\Pi(\Pi V_{q}\PP_{q})-(\Pi V_{q})\T\Pi A\Pi(\Pi V_{q})\\
  &=(V_{q}\Delta_{q}\PP_{q}^{1+\zeta})\T \Pi A\Pi( V_{q}\Delta_q\PP_{q}^{1+\zeta})-(V_{q}\Delta_{q}\PP^\zeta)\T\Pi A\Pi(V_{q}\Delta_{q}\PP^\zeta)\\
  &=(V_{q}\PP_{q}^{-1-\zeta}\Delta_{q})\T \Pi A\Pi( V_{q}\PP_{q}^{-1-\zeta}\Delta_q)-(V_{q}\PP^{-\zeta}\Delta_{q})\T\Pi A\Pi(V_{q}\PP^{-\zeta}\Delta_{q}).
\end{align*}
We remind that $\tilde E_q$ and $\PP_q$ commute and that $\Delta_q\tilde E_q=-\tilde E_q\Delta_q$ and a similar procedure for the right-hand side of \eqref{Sylvg} gives
\begin{align*}
&(\Pi V_{q}\PP_{q})\T \Pi K\Pi(\Pi V_{q}\PP_{q}\tilde E_{q})+(\Pi V_{q}\tilde E_{q})\T \Pi K\Pi(\Pi V_{q})\\
&= (V_{q}\Delta_{q}\PP_{q}^{1+\zeta})\T \Pi K\Pi(V_{q}\Delta_{q}\PP_{q}^{1+\zeta}\tilde E_{q})+( V_{q}\Delta_{q}\tilde E_{q}\PP^\zeta)\T \Pi K\Pi(V_{q}\Delta_{q}\PP_{q}^\zeta)\\
&=-(V_{q}\PP_{q}^{-1-\zeta}\Delta_{q})\T \Pi K\Pi(V_{q}\PP_{q}^{-1-\zeta} \tilde E_{q}\Delta_{q})-(V_{q}\tilde E_{q}\PP^{-\zeta}\Delta_{q})\T \Pi K\Pi(V_{q}\PP_{q}^{-\zeta}\Delta_{q}).
\end{align*}
Transposition and negation of both equations and the congruence multiplication $(\cdots)\to(\Delta_q\PP^{1+\zeta})\T(\cdots)\Delta_q\PP^{1+\zeta}$ show that the pair $(\Pi A\T \Pi,\Pi K\T\Pi)$ solves \eqref{Sylvg}, too.
Since \eqref{Sylvg} is linear, the sum or arithmetic mean of both solutions is a solution again with the persymmetric matrices $\frac12(A+\Pi A\T\Pi)$, $\frac12(K+\Pi K\T\Pi)$.
\qed
\par
For $q_1=q_2=s$ a short version of the proof applied to \eqref{Sylvs} using \eqref{PiThPi} is
\begin{align*}
 0=&\Theta\T\Pi(\Theta\T A\Theta-A-\Theta\T K\Theta E-E\T K)\Pi\Theta\\
  =&\Theta\T\big(\Theta^{-T}(\Pi A\Pi)\Theta^{-1}-(\Pi A\Pi)+\Theta^{-T}(\Pi K\Pi) E \Theta^{-1}+E\T (\Pi K\Pi)\big)\Theta.
\end{align*}
The negative transpose of the second line shows that \eqref{Sylvs} holds for $(\Pi A\T\Pi,\Pi K\T\Pi)$, as well.
\subsection{Solution structure of the Sylvester equation}
Some insight into the reasons why methods with $d_1=d_3$ may have superior properties may be gained by discussing the rank deficiencies of the Sylvester equation.
For simplicity the case $q_1=q_2=q\ge s$ is considered.
Then, the left-hand of \eqref{Sylvg} consists of a (singular) matrix mapping
\begin{align}\label{OpP}
 \PO_q:\, X\mapsto \PP_q\T X\PP_q-X,\quad X\in\R^{q\times q},
\end{align}
applied to the matrix $V_q\T A V_q$.
Now, $\PO_q$ is the mapping appearing also in the study of algebraic criteria on A-stability of Peer methods and it has been discussed in detail in \cite{Schmitt2015}.
The map $\PO_q$ is related to the maps $\LL_E:\,X\mapsto X\tilde E_q+\tilde E_q\T X$ and $\Phi_E:\, X\mapsto\int_0^1 \exp(t\tilde E_q\T)X\exp(t\tilde E_q)dt$ by $\PO_q=\Phi_E\circ \LL_E=\LL_E\circ\Phi_E$, \cite{Schmitt2015}.
Since $\Phi_E$ is a nonsingular map the kernels and images of $\PO_q$ and $\LL_E$ coincide.
\par
Since $\PO_q$ is singular, the question arises if the equation \eqref{Sylvg} for given $K$ is solvable at all.
A partial answer is given  by considering the matrix $V_q\T AV_q$ as an unknown $W$.
For $W$ the answer is affirmative since the singular factor $\tilde E_q$ appears on the right-hand side.
\begin{lemma}\label{LSylq}
For any $q\in\N$, $M\in \R^{q\times q}$ there exists a solution $W\in\R^{q\times q}$ to the equation
\begin{align}\label{PLGS}
 \PO_q(W)=\PP_q\T W\PP_q-W\stackrel!=\PP_q\T M\PP_q\tilde E_q+\tilde E_q\T M.
\end{align}
\end{lemma}
\noindent {\bf Proof:}
For stacked column vectors of $X$, the matrix associated with the map \eqref{OpP} is $\PP_q\T\otimes\PP_q\T-I_{q^2}$.
To its transpose corresponds the map $U\mapsto \PP_q U\PP_q\T-U=:\PO_q\T(U)$.
With $R:=\PP_q\T M\PP_q\tilde E_q+\tilde E_q\T M$ and by the Fredholm alternative, \eqref{PLGS} is solvable iff $tr(U\T R)=0$ for any $U$ from the kernel of $\PO_q\T$, i.e. $\PP_q U\PP_q\T=U$.
We remind that such $U$ is also in the kernel of $\LL_E\T:\,U\mapsto \tilde E_qU+U\tilde E_q\T$. Hence, with $\PO_q\T(U)=0$ it holds
\begin{align*}
 tr(U\T R)=&tr(U\T \PP_q\T M\tilde E_q\PP_q+U\T \tilde E_q\T M)\\
 =&tr\big((\PP_q U \PP_q\T)\T M\tilde E_q+U\T \tilde E_q\T M\big)\\
 =&tr\big((U\tilde E_q\T+\tilde E_q U)\T M\big)
 =tr\big(M\T \LL_E\T(U)\big)=0.
 \qed
\end{align*}
The Lemma shows that if \eqref{Sylvg} has no solution then the reason is not the singularity of $\PO_q$ but this is due to the structural restrictions on $W=V_q\T AV_q$ by the rank deficit of $V_q$, $q>s$, or the triangular form of $A$.
\par
Since solutions exist, the next question is about the solution set.
Matrices $X$ belonging to the kernel of $\PO_q$ also satisfy $\LL_E(X)=0$, which is given by the simple relations
\begin{align*}
 (j-1)x_{i,j-1}+(i-1)x_{i-1,j}=&0,\ 1\le i,j\le q,
\end{align*}
where elements with an index zero are missing.
Hence, the first $s-1$ anti-diagonals of $X$ are zero and each of the remaining anti-diagonals introduces one independent element of the kernel of $\LL_E$.
So, for the cases of most interest here, $q=3,4$, the kernels of $\LL_E$ are given by
\begin{align}\label{SylKern}
 X_3=\begin{pmatrix}
  0&0&2\xi_1\\
  0&-\xi_1&\xi_2\\ 2\xi_1&-\xi_2&\xi_3
 \end{pmatrix},\quad
 X_4=\begin{pmatrix}
  0&0&0&3\xi_1\\ 0&0&-\xi_1&3\xi_2\\ 0&\xi_1&-2\xi_2&\xi_3\\-3\xi_1&3\xi_2&-\xi_3&\xi_4
 \end{pmatrix}\,.
\end{align}
For $q=s$ the matrix $V_q$ is nonsingular and the kernel of equation
\eqref{Sylvg} is easily found.
For practical methods, lower triangular form is of interest.
\begin{lemma}\label{LKern3}
For $q=s=3$ the matrix map $X\mapsto \PO_3(V_3\T XV_3)$ on the left-hand side of \eqref{Sylvg} has non-trivial kernel elements with matrices in lower triangular form iff $d_1=d_3$.
In this case the kernel is spanned by the single matrix
\[ \begin{pmatrix} 1&0&0\\ -3&1&0\\3&-3&1 \end{pmatrix}.
\]
\end{lemma}
\noindent {\bf Proof:}
With $X_3$ from \eqref{SylKern}, the kernel of $\PO_3(V_3\T XV_3)$ is given by $\hat X:=V_3^{-T}X_3V_3^{-1}$.
The conditions that all super-diagonals vanish are given by the linear system
\[ \begin{pmatrix}
  (d_1+d_3)d_3&d_1&-1\\
  d_1(d_1+d_3)&d_3&-1\\
  -d_1d_3&d_1+d_3&-1
\end{pmatrix}\begin{pmatrix} \xi_1\\\xi_2\\\xi_3\end{pmatrix}=0.
\]
Nontrivial solutions with $d_j\not=0$ exist only if the determinant $d_1d_3(d_3-d_1)$ vanishes with $d_3=d_1$.
Then, up to factors kernel elements are multiples of the matrix from the statement.
\qed
\par
\begin{remark}
This lemma may give a first hint why methods with equal node differences obtain higher orders since the loss of one degree of freedom by fixing the parameters $d_3=d_1$ is compensated for by the free factor of the kernel element for $q=3$.
This happens only once, for order $q=4=s+1$ the kernel is trivial.
\end{remark}

\section{Numerical Results}\label{STests}
We present numerical results for three different methods:
\begin{center}
\begin{tabular}{ll}
\hline
\rule{0mm}{0.4cm}Name & coefficients \\
\hline
\rule{0mm}{0.4cm}\verb+BDF3o22+ & (\ref{BDF3n})--(\ref{BDF3s}) \\
\verb+BDF3o32+ & (\ref{BDF3n}), (\ref{ANimp})--(\ref{A0imp}) \\
\verb+PEER3o32w+ & (\ref{BDF3n}), (\ref{BDF3o32N})--(\ref{BDF3o32s}) \\
\hline
\end{tabular}
\end{center}
\par
All calculations have been done with Matlab-Version R2019a, using the
nonlinear solver \textit{fsolve} to approximate the overall coupled scheme
(\ref{ADJEnd})--(\ref{FORWEnd}) with a tolerance $1e\!-\!14$.
To illustrate the rates of convergence, we consider two unconstrained nonlinear
optimal control problems.

\subsection{The Rayleigh problem}
The first problem is taken from \cite{JacobsonMayne1970} and describes the behaviour of
a tunnel-diode oscillator. With the electric current $y_1(t)$ and the
transformed voltage at the generator $u(t)$, the unconstrained Rayleigh
problem reads
\begin{align}
\mbox{Minimize } \int_0^{2.5} u(t)^2+y_1(t)^2\,dt \label{rayleigh_objfunc} &\\
\mbox{subject to } y''_1(t)-y'_1\left( 1.4 - 0.14 y'_1(t)^2 \right)+y_1(t)
=&\,  4u(t),\quad t\in(0,2.5], \label{rayleigh_ODE}\\
y_1(0)=y'_1(0)=&\,-\!5. \label{rayleigh_ODEinit}
\end{align}
Introducing $y_2(t)=y_1'(t)$ and eliminating the control $u(t)$ yields
the following nonlinear boundary value problem (see \cite{LangVerwer2013}
for more details):
\begin{align}
y'_1(t) =&\, y_2(t),\label{rayleigh-bvp1}\\
y'_2(t) =&\, -y_1(t) + y_2\left( 1.4 - 0.14 y_2(t)^2 \right) - 8p_2(t),\\
         &\, y_1(0) = -5,\,y_2(0)=-5,\\
p'_1(t) =&\, p_2(t) - 2 y_1(t),\\
p'_2(t) =&\, -p_1(t) - (1.4-0.42 y_2(t)^2)p_2(t),\\
            &\, p_1(2.5) = 0,\,p_2(2.5)=0.\label{rayleigh-bvp2}
\end{align}
To study convergence orders of our new methods, we compute a reference solution
by applying the classical fourth-order RK4 with $N\!=\!320$. Numerical results
are presented in Table~\ref{tab-rayleigh}. The two BDF3 methods give the same
results and perform by nearly a factor of three better than the PEER3 method
(having $c_3<1$) in terms of errors.
As expected and supported by the theory, the convergence orders for the state
variables are
nearly three and range between two and three for the adjoint variables, see Remark~\ref{RKonvp}.
Surprisingly, even method \verb+BDF3o22+ which misses the highest order condition in the end step draws level with \verb+BDF3o32+.

\begin{table}[ht!]
\begin{center}
\begin{tabular}{|l||r|r|r|r|r|r|}\hline
\rule{0mm}{0.4cm}N & 40 & 80 & 160 & 320 \\ \hline \hline%
\rule{0mm}{0.4cm} BDF3o22 &&&&\\
\rule{0mm}{0.4cm} $y_1$ (order)
&$4.23e\!-\!4$&$5.67e\!-\!5$ (2.9)&$7.68e\!-\!6$ (2.9)&$8.98e\!-\!7$ (3.1)\\
\rule{0mm}{0.4cm} $y_2$ (order)
&$7.05e\!-\!3$&$1.39e\!-\!3$ (2.3)&$2.19e\!-\!4$ (2.7)&$3.08e\!-\!5$ (2.8)\\
\rule{0mm}{0.4cm} $p_1$ (order)
&$1.65e\!-\!3$&$2.63e\!-\!4$ (2.6)&$4.76e\!-\!5$ (2.5)&$9.16e\!-\!6$ (2.4)\\
\rule{0mm}{0.4cm} $p_2$ (order)
&$3.45e\!-\!2$&$6.79e\!-\!3$ (2.3)&$1.58e\!-\!3$ (2.1)&$3.89e\!-\!4$ (2.0)\\
\hline
\rule{0mm}{0.4cm} BDF3o32 &&&&\\
\rule{0mm}{0.4cm} $y_1$ (order)
&$4.23e\!-\!4$&$5.67e\!-\!5$ (2.9)&$7.68e\!-\!6$ (2.9)&$8.98e\!-\!7$ (3.1)\\
\rule{0mm}{0.4cm} $y_2$ (order)
&$7.05e\!-\!3$&$1.39e\!-\!3$ (2.3)&$2.19e\!-\!4$ (2.7)&$3.08e\!-\!5$ (2.8)\\
\rule{0mm}{0.4cm} $p_1$ (order)
&$1.65e\!-\!3$&$2.63e\!-\!4$ (2.6)&$4.76e\!-\!5$ (2.5)&$9.16e\!-\!6$ (2.4)\\
\rule{0mm}{0.4cm} $p_2$ (order)
&$3.45e\!-\!2$&$6.79e\!-\!3$ (2.3)&$1.58e\!-\!3$ (2.1)&$3.89e\!-\!4$ (2.0)\\
\hline
\rule{0mm}{0.4cm} PEER3o32w &&&&\\
\rule{0mm}{0.4cm} $y_1$ (order)
&$1.75e\!-\!3$&$2.13e\!-\!4$ (3.0)&$2.60e\!-\!5$ (3.0)&$2.99e\!-\!6$ (3.1)\\
\rule{0mm}{0.4cm} $y_2$ (order)
&$6.01e\!-\!3$&$8.96e\!-\!4$ (2.8)&$1.22e\!-\!4$ (2.9)&$1.53e\!-\!5$ (3.0)\\
\rule{0mm}{0.4cm} $p_1$ (order)
&$3.75e\!-\!3$&$6.12e\!-\!4$ (2.6)&$1.30e\!-\!4$ (2.2)&$2.92e\!-\!5$ (2.2)\\
\rule{0mm}{0.4cm} $p_2$ (order)
&$9.96e\!-\!2$&$2.45e\!-\!2$ (2.0)&$5.92e\!-\!3$ (2.0)&$1.45e\!-\!3$ (2.0)\\
\hline
\end{tabular}
\vspace{-2mm}
\parbox{13cm}{
\caption{Rayleigh problem: $l^\infty$-convergence of the discrete state
errors $y_i(t_n)-Y_{ni}$ and adjoint state errors $p_i(t_n)-P_{ni}$ for
\texttt{BDF3o22}, \texttt{BDF3o32}, and \texttt{PEER3o32w}.
The numbers in brackets estimate the order of convergence.
\label{tab-rayleigh}}}
\end{center}
\end{table}

\subsection{The van der Pol problem}
\begin{table}[ht!]
\begin{center}
\begin{tabular}{|l||r|r|r|r|r|}\hline
\rule{0mm}{0.4cm}N & 160 & 320 & 640 & 1280 \\ \hline \hline%
\rule{0mm}{0.4cm} BDF3o22 &&&&\\
\rule{0mm}{0.4cm} $x_1$ (order)
&$1.01e\!-\!5$&$1.34e\!-\!6$ (2.9)&$1.73e\!-\!7$ (3.0)&$2.39e\!-\!8$ (2.9)\\
\rule{0mm}{0.4cm} $x_2$ (order)
&$8.26e\!-\!6$&$1.07e\!-\!6$ (3.0)&$1.39e\!-\!7$ (2.9)&$1.77e\!-\!8$ (3.0)\\
\rule{0mm}{0.4cm} $p_1$ (order)
&$7.92e\!-\!3$&$1.91e\!-\!3$ (2.0)&$4.68e\!-\!4$ (2.0)&$1.16e\!-\!4$ (2.0)\\
\rule{0mm}{0.4cm} $p_2$ (order)
&$7.32e\!-\!3$&$1.77e\!-\!3$ (2.0)&$4.32e\!-\!4$ (2.0)&$1.07e\!-\!4$ (2.0)\\
\hline
\rule{0mm}{0.4cm} BDF3o32 &&&&\\
\rule{0mm}{0.4cm} $x_1$ (order)
&$1.01e\!-\!5$&$1.34e\!-\!6$ (2.9)&$1.73e\!-\!7$ (3.0)&$2.39e\!-\!8$ (2.9)\\
\rule{0mm}{0.4cm} $x_2$ (order)
&$8.26e\!-\!6$&$1.07e\!-\!6$ (3.0)&$1.39e\!-\!7$ (2.9)&$1.77e\!-\!8$ (3.0)\\
\rule{0mm}{0.4cm} $p_1$ (order)
&$7.92e\!-\!3$&$1.91e\!-\!3$ (2.0)&$4.68e\!-\!4$ (2.0)&$1.16e\!-\!4$ (2.0)\\
\rule{0mm}{0.4cm} $p_2$ (order)
&$7.32e\!-\!3$&$1.77e\!-\!3$ (2.0)&$4.32e\!-\!4$ (2.0)&$1.07e\!-\!4$ (2.0)\\
\hline
\rule{0mm}{0.4cm} PEER3o32w &&&&\\
\rule{0mm}{0.4cm} $x_1$ (order)
&$2.19e\!-\!5$&$3.25e\!-\!6$ (2.8)&$4.42e\!-\!7$ (2.9)&$6.21e\!-\!8$ (2.8)\\
\rule{0mm}{0.4cm} $x_2$ (order)
&$9.76e\!-\!6$&$1.23e\!-\!6$ (3.0)&$1.54e\!-\!7$ (3.0)&$1.94e\!-\!8$ (3.0)\\
\rule{0mm}{0.4cm} $p_1$ (order)
&$2.42e\!-\!2$&$6.35e\!-\!3$ (1.9)&$1.62e\!-\!3$ (2.0)&$4.11e\!-\!4$ (2.0)\\
\rule{0mm}{0.4cm} $p_2$ (order)
&$2.24e\!-\!2$&$5.86e\!-\!3$ (1.9)&$1.50e\!-\!3$ (2.0)&$3.80e\!-\!4$ (2.0)\\
\hline
\end{tabular}
\vspace{-2mm}
\parbox{13cm}{
\caption{Van der Pol problem with $\varepsilon=0.1$: $l^\infty$-convergence of the discrete state errors $y_i(t_n)-Y_{ni}$ and adjoint state errors $p_i(t_n)-P_{ni}$
for \texttt{BDF3o22}, \texttt{BDF3o32}, and \texttt{PEER3o32w}.
The numbers in brackets estimate the order of convergence.}
\label{tab-vdpol}}
\end{center}
\end{table}
The second example is the following optimal control problem for the
van der Pol oscillator:
\begin{align}
\mbox{Minimize } \int_0^2 u(t)^2+y(t)^2+y'(t)^2\,dt \label{vanderpol_objfunc} &\\
\mbox{subject to } \;\varepsilon y''(t)-(1-y(t)^2)y'(t)+y(t)
=&\, u(t),\quad t\in(0,2], \label{vanderpol_ODE}\\
y(0)=0,\;y'(0)=&\, 2. \label{vanderpol_ODEinit}
\end{align}
We set $\varepsilon=0.1$ and use Lienhard's coordinates
$y_2(t)=y(t)$, $y_1(t)=\varepsilon y'(t)+y(t)^3/3-y(t)$ to
end up with the boundary value problem (see \cite{LangVerwer2013}
for more details)
\begin{align}
y'_1(t) =&\, -y_2(t)-\frac{p_1(t)}{2},\label{vanderpol-bvp1}\\
y'_2(t) =&\, \frac{1}{\varepsilon}
\left( y_1(t)+y_2(t)-\frac{y_2(t)^3}{3}\right),\\
&\, y_1(0) = 2\varepsilon, \,y_2(0)=0,\\
p'_1(t) =&\, -\frac{1}{\varepsilon}p_2(t) -
\frac{2}{\varepsilon^2} \left( y_1(t)+y_2(t)-\frac{y_2(t)^3}{3}\right),\\
p'_2(t) =&\, p_1(t)-\frac{1}{\varepsilon}\left( 1-y_2(t)^2\right)p_2(t)\nonumber\\
&\, -\frac{2}{\varepsilon^2} \left( y_1(t)+y_2(t)-\frac{y_2(t)^3}{3}\right)
\left( 1-y_2(t)^2\right)-2y_2(t),\\
&\, p_1(2) = 0,\,p_2(2)=0.\label{vanderpol-bvp2}
\end{align}
For comparison, a reference solution is computed for $N=2560$ using the
W-method ROS3WO from \cite{LangVerwer2013}. In Table~\ref{tab-vdpol}, numerical
results for $N=160,320,640$, and $1280$ are shown.
Obviously, the two BDF3 methods once again deliver equal results and outperform the PEER method by a factor three in terms of errors despite the lower order in the end step of \verb+BDF3o22+.
Here, the convergence
orders three for the state variables and two for the adjoint variables according to Theorem~\ref{TKonv} are visible quite clearly.

\section{Summary}
By introducing a redundant formulation of Peer two-step methods and exceptional boundary steps, sufficient additional degrees of freedom could be gained to prove order $s\!=\!3$ for the state solution and $s\!-\!1\!=\!2$ for the adjoint variables of the full boundary value problem derived from the first-order optimality conditions.
Although a detailed analysis for the global order pair (3,2) detected some exotic schemes being A-stable, the most attractive standard Peer methods in the interior of the grid are based on the BDF3 scheme.
Different approaches for the adjoint boundary condition lead to three methods which reproduce the correct orders in numerical tests with two nonlinear problems.
Some matrix background helps to explain why flip symmetry of the nodes of BDF may lead to its superior properties here.

\vspace{0.5cm}
\par
\noindent {\bf Acknowledgements.}
The first author is supported by the Deutsche Forschungsgemeinschaft
(DFG, German Research Foundation) within the collaborative research center
TRR154 {\em ``Mathematical modeling, simulation and optimisation using
the example of gas networks''} (Project-ID 239904186, TRR154/2-2018, TP B01).

\bibliographystyle{plain}
\bibliography{bibpeeropt}

\end{document}